\documentclass[final,leqno,onefignum,onetabnum]{siamltex1213}
\usepackage{amsmath,amssymb,boxedminipage,algorithm,algorithmic}
\usepackage{epsfig,graphicx,epstopdf}

\title{Some Results on the Regularization of LSQR for Large-Scale Discrete
Ill-Posed Problems\thanks{This work was supported in part by the
National Basic Research Program
of China 2011CB302400 and the National Science Foundation of China (No. 11371219)}}

\author{Yi Huang\thanks{Department of Mathematical Sciences, Tsinghua
       University, Beijing 100084, People's Republic of China.
       (huangyi10@mails.tsinghua.edu.cn)}
       \and Zhongxiao Jia\thanks{Corresponding author. Department
       of Mathematical Sciences, Tsinghua
        University, Beijing 100084, People's Republic of China.
        (jiazx@tsinghua.edu.cn)}}

\begin{document}
\maketitle
\slugger{simax}{xxxx}{xx}{x}{x--x}

\begin{abstract}
LSQR, a Lanczos bidiagonalization based Krylov subspace iterative
method, and its mathematically equivalent CGLS applied to normal equations
system, are commonly used for large-scale discrete ill-posed problems.
It is well known that LSQR and CGLS have
regularizing effects, where the number of iterations plays the role
of the regularization parameter. However, it has long been unknown whether the
regularizing effects are good enough to find best possible regularized
solutions. Here a best possible regularized solution means that it is at least as
accurate as the best regularized solution obtained by the truncated singular
value decomposition (TSVD) method. In this paper, we establish
bounds for the distance between the $k$-dimensional Krylov subspace and
the $k$-dimensional dominant right singular space.
They show that the Krylov subspace captures
the dominant right singular space better for severely and moderately
ill-posed problems than for mildly ill-posed problems. Our general conclusions
are that LSQR has better regularizing effects for the first two
kinds of problems than for the third kind, and a hybrid LSQR with
additional regularization is generally needed for mildly ill-posed problems.
Exploiting the established bounds, we derive an estimate for the accuracy of
the rank $k$ approximation generated by Lanczos
bidiagonalization. Numerical experiments illustrate that the regularizing effects
of LSQR are good enough to compute best possible regularized solutions
for severely and moderately ill-posed problems, stronger than our theory
predicts, but they are not for mildly
ill-posed problems and additional regularization is needed.
\end{abstract}



\pagestyle{myheadings}
\thispagestyle{plain}
\markboth{YI HUANG AND ZHONGXIAO JIA}{SOME RESULTS ON REGULARIZATION OF LSQR}

\section{Introduction}
We consider the iterative solution of large-scale discrete
ill-posed problems
\begin{equation}\label{eq1}
  \min\limits_{x\in \mathbb{R}^{n}}\|Ax-b\|,\ \ \ A\in \mathbb{R}^{m\times n},
  \ b\in \mathbb{R}^{m},\ \ m\geq n,
\end{equation}
where the norm $\|\cdot\|$ is the 2-norm of a vector or matrix,
and the matrix $A$ is extremely ill conditioned with its singular values decaying
gradually to zero without a noticeable gap. This kind of problem arises in many
science and engineering areas, such as signal processing and image restoration,
typically when discretizing Fredholm integral equations of the first-kind
\cite{hansen98,hansen10}. In particular, the right-hand side
$b$ is affected by noise, caused by measurement or discretization errors, i.e.,
\begin{equation*}\label{}
  b=\hat{b}+e,
\end{equation*}
where $e\in \mathbb{R}^{m}$ represents the Gaussian white noise vector and
$\hat{b}\in \mathbb{R}^{m}$ denotes the noise-free right-hand side, and it
is supposed that $\|e\|<\|\hat{b}\|$.
Because of the presence of noise $e$ in $b$ and the ill-conditioning of $A$,
the naive solution $x_{naive}=A^{\dagger}b$ of \eqref{eq1} is meaningless and
far from the true solution $x_{true}=A^{\dagger}\hat{b}$, where the superscript
$\dagger$ denotes the Moore-Penrose generalized inverse of a matrix.
Therefore, it is necessary to use regularization to determine a
best possible approximation to $x_{true}=A^{\dagger}\hat{b}$ \cite{engl96,hanke95,
hansen98,hansen10}.

The solution of \eqref{eq1} can be analyzed by the SVD of $A$:
\begin{equation}\label{eqsvd}
  A=U\left(\begin{array}{c} \Sigma \\ \mathbf{0} \end{array}\right) V^{T},
\end{equation}
where $U = (u_1,u_2,\ldots,u_m)\in\mathbb{R}^{m\times m}$ and
$V = (v_1,v_2,\ldots,v_n)\in\mathbb{R}^{n\times n}$ are orthogonal matrices,
and the entries of the diagonal matrix $\Sigma = \mathrm{diag}(\sigma_1,\sigma_2,
\ldots,\sigma_n)\in\mathbb{R}^{n\times n}$ are the singular values of $A$, which
are assumed to be simple throughout the paper and labelled in decreasing order
$\sigma_1>\sigma_2 >\cdots >\sigma_n>0$. With \eqref{eqsvd}, we obtain
\begin{equation}\label{eq4}
  x_{naive}=\sum\limits_{i=1}^{n}\frac{u_i^{T}b}{\sigma_i}v_i =
  \sum\limits_{i=1}^{n}\frac{u_i^{T}\hat{b}}{\sigma_i}v_i +
  \sum\limits_{i=1}^{n}\frac{u_i^{T}e}{\sigma_i}v_i
  =x_{true}+\sum\limits_{i=1}^{n}\frac{u_i^{T}e}{\sigma_i}v_i.
\end{equation}
Throughout the paper, we assume that $\hat{b}$ satisfies the discrete Picard
condition: On average, the coefficients $\mid u_i^{T}\hat{b}\mid$
decay faster than the singular values. To be definitive, for simplicity
we assume that these coefficients satisfy a widely used model in the literature,
e.g., \cite[p. 81, 111 and 153]{hansen98} and \cite[p. 68]{hansen10}:
\begin{equation}\label{picard}
  \mid u_i^T \hat b\mid=\sigma_i^{1+\beta},\ \ \beta>0,\ i=1,2,\ldots,n.
\end{equation}

Let $k_0$ be the transition point such that $|u_{k_0}^T \hat b|>
|u_{k_0+1}^T e |$ and $| u_{k_0+1}^T \hat b |\leq
|u_{k_0+1}^T e|$ \cite[p. 98]{hansen10}. Then the TSVD method computes
\begin{equation*}\label{}
  x^{TSVD}_k=\left\{\begin{array}{ll} \sum\limits_{i=1}^{k}\frac{u_i^{T}b}
  {\sigma_i}{v_i}\thickapprox
  \sum\limits_{i=1}^{k}\frac{u_i^{T}\hat{b}}
{\sigma_i}{v_i},\ \ \ &k\leq k_0;\\ \sum\limits_{i=1}^{k}\frac{u_i^{T}b}
{\sigma_i}{v_i}\thickapprox
\sum\limits_{i=1}^{k_0}\frac{u_i^{T}\hat{b}}{\sigma_i}{v_i}+
\sum\limits_{i=k_0+1}^{k}\frac{u_i^{T}e}{\sigma_i}{v_i},\ \ \ &k>k_0,
\end{array}\right.
\end{equation*}
which can be written as $x_k^{TSVD}=A_k^{\dagger}b$, the solution of the
modified problem that replaces $A$ by its best rank $k$ approximation
$A_k=U_k\Sigma_k V_k^T$ in (\ref{eq1}), where $U_k=(u_1,\ldots,u_k)$,
$V_k=(v_1,\ldots,v_k)$ and $\Sigma_k={\rm diag}(\sigma_1,\ldots,\sigma_k)$.
Remarkably, $x_{k}^{TSVD}$ is the minimum-norm least squares solution of
the perturbed problem that replaces $A$ in \eqref{eq1} by its best rank $k$
approximation $A_{k}$, and the best possible TSVD solution
of \eqref{eq1} by the TSVD method is $x_{k_0}^{TSVD}$ \cite[p. 98]{hansen10}.
A number of approaches have been proposed for determining $k_0$, such as discrepancy
principle, discrete L-curve and generalized cross validation; see, e.g.,
\cite{bauer11,bazan10,hansen98,kilmer03,reichel13}
for comparisons of the classical and new ones. In our numerical experiments, we
use the L-curve criterion in the TSVD method and hybrid LSQR. The TSVD method has
been widely studied; see, e.g.,
\cite{bergou14,hansen98,hansen10,li11}.

For a small and moderate \eqref{eq1}, the TSVD method has been used as
a general-purpose reliable and efficient numerical
method for solving \eqref{eq1}. As a result, we will take the TSVD
solution $x_{k_0}^{TSVD}$ as a standard reference when assessing the regularizing
effects of iterative
solvers and accuracy of iterates under consideration in this paper.

As well known, it is generally not feasible to compute SVD when \eqref{eq1}
is large. In this case, one typically projects \eqref{eq1} onto a
sequence of low dimensional Krylov subspaces and gets a sequence of iterative
solutions \cite{hanke95,hansen98,hansen10,vogel02}.
The Conjugate Gradient (CG) method has been used when $A$ is symmetric
definite \cite{hanke95}. As a CG-type method applied to the semidefinite
linear system $Ax=b$ or the normal equations system $A^TAx=A^Tb$,
the CGLS algorithm has been studied; see \cite{bjorck96,hansen98,hansen10}
and the references therein. The LSQR algorithm \cite{paige82}, which
is mathematically equivalent to CGLS, has attracted great attention, and
is known to have regularizing effects and exhibits semi-convergence
(see \cite[p. 135]{hansen98}, \cite[p. 110]{hansen10}, and also
\cite{bjorck88,hanke01,hnetyn09,oleary81}): The iterates
tend to be better and better approximations to the exact solution $x_{true}$
and their norms increase slowly and the residual norms decrease.
In later stages, however, the noise $e$ starts to deteriorate the iterates,
so that they will start to diverge from
$x_{true}$ and instead converge to the naive solution $x_{naive}$,
while their norms increase considerably and the residual norms stabilize.
Such phenomenon is due to the fact that a
projected problem inherits the ill-conditioning of \eqref{eq1}.
That is, as the iterations proceed, the noise progressively enters
the solution subspace, so that a small singular value of the projected
problem appears and the regularized solution is deteriorated.

As far as an iterative solver for solving \eqref{eq1} is concerned,
a central problem is whether or not a pure iterative solver has already obtained
a best possible regularized solution at semi-convergence, namely
whether or not the regularized solution at semi-convergence is at
least as accurate as $x_{k_0}^{TSVD}$.
As it appears, for Krylov subspace based iterative solvers, their
regularizing effects
critically rely on how well the underlying $k$-dimensional Krylov subspace
captures the $k$-dimensional dominant right singular subspace of $A$.
The richer information the Krylov subspace contains on the $k$-dimensional dominant
right singular subspace,
the less possible a small Ritz value of the resulting projected problem appears
and thus the better regularizing effects the solver has. To precisely describe
the regularizing effects of an iterative solver, we introduce the term of
{\em full} or {\em partial}
regularization: If the iterative solver itself computes a best possible
regularized solution at semi-convergence, it is said
to have the full regularization; in this case, no additional regularization
is needed. Here, as defined in the abstract,
a best possible regularized solution means that it is at least as
accurate as the best regularized solution obtained by the truncated singular
value decomposition (TSVD) method.
Otherwise, it is said to have the partial regularization;
in this case, in order to compute a best possible regularized solution,
its hybrid variant, e.g., a hybrid LSQR, is needed that combines the solver
with additional regularization
\cite{bjorck88,chung08,lewis09,neuman12,novati13,oleary81},
which aims to remove the effects of small Ritz values,
and expand the Krylov subspace until it captures all the dominant
SVD components needed and the method obtains
a best possible regularized solution.
The study of the regularizing effects of LSQR
and CGLS has been receiving intensive attention for years;
see \cite{hansen98,hansen10} and the references therein. However,
there has yet been no definitive result or assertion on their full or partial
regularization.

To proceed, we need
the following definition of the degree of ill-posedness, which
follows Hofmann's book \cite{hofmann86} and has been commonly used
in the literature, e.g.,
\cite{hansen98,hansen10}: If there exists a positive real number
$\alpha$ such that the singular values satisfy $\sigma_j=
\mathcal{O}(j^{-\alpha})$, then the problem is
termed as mildly or moderately ill-posed if $\alpha\le1$ or $\alpha>1$;
if $\sigma_j=\mathcal{O}(e^{-\alpha j})$ with $\alpha>0$ considerably,
$j=1,2,\ldots,n$, then the problem is termed severely ill-posed.
It is clear that the singular values $\sigma_j$ of a
severely ill-posed problem decay exponentially at the same rate $\rho^{-1}$,
while those of a moderately or mildly ill-posed problem decay more and more
slowly at the decreasing rate $\left(\frac{j}{j+1}\right)^{\alpha}$
approaching one with increasing $j$, which, for the same $j$, is
smaller for the moderately ill-posed problem than it
for the mildly ill-posed problem.

Other minimum-residual methods have also gained attention for solving \eqref{eq1}.
For problems with $A$ symmetric, MINRES and its preferred variant
MR-II are alternatives and
have been shown to have regularizing effects \cite{hanke95}. When $A$ is
nonsymmetric and multiplication with $A^{T}$ is difficult or impractical to
compute, GMRES and its preferred variant RRGMRES are candidates
\cite{calvetti02_1,neuman12}. The hybrid approach
based on the Arnoldi process was first introduced in \cite{calvetti00},
and has been studied in
\cite{calvetti01,calvetti00,calvetti03,lewis09}. Recently, Gazzola {\em et al}.
\cite{gazzola13,gazzola14,gazzola-online,novati13}
have studied more methods based on the Lanczos bidiagonalization, the Arnoldi
process and the nonsymmetric Lanczos process for the severely ill-posed
problem \eqref{eq1}. They have described a general framework of the hybrid
methods and present Krylov-Tikhonov methods with different parameter choice
strategies employed.

In this paper, we focus on LSQR. We derive bounds for the 2-norm distance between
the underlying $k$-dimensional Krylov subspace and the $k$-dimensional right
singular space. There has been no rigorous and quantitative result on
the distance before. The results indicate that the $k$-dimensional Krylov
subspace captures the $k$-dimensional dominant right singular space better
for severely and moderately ill-posed problems than for mildly ill-posed problems.
As a result, LSQR has better regularizing effects for the first two kinds of
problems than for the third kind. By the bounds and the
analysis on them, we draw a definitive conclusion that LSQR generally has only
the partial regularization for mildly ill-posed problems, so that a hybrid LSQR with
additional explicit regularization is needed to compute a best possible regularized
solution. We also use the bounds to derive an estimate for the accuracy of the
rank $k$ approximation, generated by Lanczos bidiagonalization, to $A$, which
is closely related to the regularization of LSQR.
Our results help to further understand the regularization of LSQR,
though they appear less sharp. In addition, we derive a bound on the diagonal entries
of the bidiagonal matrices generated by the Lanczos bidigonalization process,
showing how fast they decay.
Numerical experiments confirm our theory that LSQR has only the partial
regularization for mildly ill-posed problems and a hybrid LSQR is needed
to compute best possible regularized solutions. Strikingly,
the experiments demonstrate that LSQR
has the full regularization for severely and moderately ill-posed problems.
Our theory gives a partial support for the observed general
phenomena. Throughout the paper,
all the computation is assumed in exact arithmetic. Since CGLS is
mathematically equivalent to LSQR, all the assertions on LSQR apply to CGLS.

This paper is organized as follows. In Section \ref{SectionMain},
we describe the LSQR algorithm, and then present our theoretical results
on LSQR with a detailed analysis. In Section \ref{SectionExp}, we report
numerical experiments to justify the partial regularization of LSQR
for mildly ill-posed problems. We also report some definitive and general
phenomena observed. Finally, we conclude the paper in Section \ref{SectionCon}.

Throughout the paper, we denote by $\mathcal{K}_{k}(C, w)=
span\{w,Cw,\ldots,C^{k-1}w\}$ the $k$-dimensional Krylov subspace generated
by the matrix $\mathit{C}$ and the vector $\mathit{w}$, by $\|\cdot\|_F$
the Frobenius norm of a matrix, and by $I$ the identity matrix with order clear
from the context.

\section{The regularization of LSQR}\label{SectionMain}

LSQR for solving \eqref{eq1} is based
on the Lanczos bidiagonalization process, which starts with
$p_1=b/\|b\|$ and, at step (iteration) $k$,
computes two orthonormal
bases $ \{q_1,q_2,\dots,q_k\}$ and $ \{p_1,p_2,\dots,p_k\}$ of the Krylov
subspaces $\mathcal{K}_{k}(A^{T}A,A^{T}b)$ and $\mathcal{K}_{k}(A A^{T},b)$,
respectively.


Define the matrices $Q_k=(q_1, q_2, \ldots, q_k)$ and $P_{k+1}=(p_1, p_2, \ldots, p_{k+1})$.
Then the $k$-step Lanczos bidiagonalization can be written in the matrix form
\begin{eqnarray}
  AQ_k&=&P_{k+1}B_k, \label{eqmform1}\\
  A^TP_{k+1}&=&Q_{k}B_k^T+\alpha_{k+1}q_{k+1}e_{k+1}^{T},\label{eqmform2}
\end{eqnarray}
where $e_{k+1}$ denotes the $(k+1)$-th canonical basis vector
of $\mathbb{R}^{k+1}$ and the quantities $\alpha_i, i=1,2,\ldots,k+1$ and
$\beta_i,\ i=2,\ldots,k+1$ denote the diagonal
and subdiagonal elements of the $(k+1)\times k$ lower bidiagonal matrix $B_k$,
respectively.
At iteration $k$, LSQR computes the solution $x^{(k)}=Q_ky^{(k)}$ with
\begin{equation*}\label{}
  y^{(k)}=\arg\min\limits_{y\in \mathbb{R}^{k}}\|{\|b\|e_1 - B_ky}\|.
\end{equation*}
Note that
$P_{k+1}^Tb=\|b\|e_1$. We get
\begin{equation}\label{lsqrsol}
x^{(k)}=Q_ky^{(k)}=\|b\|Q_kB_k^{\dagger}e_1=Q_kB_k^{\dagger}P_{k+1}^Tb.
\end{equation}

As stated in the introduction, LSQR exhibits semi-convergence at some
iteration: The iterates $x^{(k)}$ become better
approximations to $x_{true}$ until some iteration $k$, and the noise
will dominate the $x^{(k)}$ after that iteration. The iteration number $k$
plays the role of the regularization parameter.
However, semi-convergence does not necessarily mean that LSQR finds a best
possible regularized solution as $B_k$ may become ill-conditioned before $k\leq k_0$
but $x^{(k)}$ does not yet contain all the needed $k_0$ dominant SVD components
of $A$. In this case, in order to get a best possible regularized solution,
one has to use a hybrid LSQR method, as described in the introduction.
The significance of \eqref{lsqrsol} is that the LSQR iterates can be
interpreted as the minimum-norm least squares solutions of the perturbed
problems that replace $A$ in \eqref{eq1} by its rank $k$ approximations
$Q_kB_k^{\dagger}P_{k+1}^T$, whose nonzero singular values are just those
of $B_k$. If the singular values of $B_k$ approximate
the $k$ large singular values of $A$ in natural order for $k=1,2,\ldots,k_0$,
then LSQR must have the full regularization, and the regularized solution
$x^{(k_0)}$ is best possible and is as comparably accurate as the best possible
regularized solution $x_{k_0}^{TSVD}$ by the TSVD method.

Hansen's analysis \cite[p. 146]{hansen98} shows that the LSQR
  iterates have the filtered SVD expansions:
  \begin{equation*}
  {x}^{(k)}=\sum\limits_{i=1}^nf_i^{(k)}\frac{u_i^{T}b}{\sigma_i}v_i,
  \end{equation*} where $f_i^{(k)}=1-\prod\limits_{j=1}^k\frac{\theta_j^{(k)}
  -\sigma_i}
  {\theta_j^{(k)}},\ i=1,2,\ldots,n$, and $\theta_j^{(k)},\
  j=1,2,\ldots,k$ are the singular values of $B_k$.
  In our context, if we have $\theta_{k}^{(k)}\leq \sigma_{k_0+1}$
  for some $k\leq k_0$, the factors
  $f_i^{(k)}$, $i=k+1,\ldots,n$ are not small, meaning that $x^{(k)}$
  is already deteriorated and becomes a poorer regularized solution, namely,
  LSQR surely does not have full regularization.
  As a matter of fact, in terms of the best possible solution $x_{k_0}^{TSVD}$,
  it is easily justified that the full regularization of LSQR
  is equivalent to requiring that the singular values of $B_k$ approximate
  the $k$ largest singular values of $A$ in natural order for $k=1,2,\ldots,k_0$,
  so it is impossible
  to have $\theta_{k}^{(k)}\leq \sigma_{k_0+1}$ for $k\leq k_0$.

The regularizing effects of LSQR
critically depend on what $\mathcal{K}_k(A^TA,A^Tb)$ mainly contains
and provides. Note that the eigenpairs of $A^TA$ are the squares of
singular values and right singular vectors of $A$, and the
tridiagonal matrix $B_k^TB_k$ is the
projected matrix of $A^TA$ onto the subspace $\mathcal{K}_k(A^TA,A^Tb)$,
which is obtained by applying the symmetric Lanczos tridiagonalization
process to $A^TA$ starting with $q_1=A^Tb/\|A^Tb\|$ \cite{bjorck96}.
We have a general claim deduced from \cite{bjorck96,parlett} and exploited
widely in \cite{hansen98,hansen10}:
The more information the subspace $\mathcal{K}_k(A^TA,A^Tb)$ contains
on the $k$ dominant right singular vectors,
the more possible and accurate the $k$ Ritz values approximate
the $k$ largest singular values of $A$; on the other hand,
the less information it contains on the other $n-k$ right singular vectors,
the less accurate a small Ritz value is if it appears.
For our problem, since the small singular values of $A$ are clustered and
close to zero, it is expected that a small Ritz value will show up as $k$
grows large, and it starts to appear more late when $\mathcal{K}_k(A^TA,A^Tb)$
contains less information on the other $n-k$ right singular vectors.
In this sense, we say that LSQR has better regularizing effects since $x^{(k)}$
contains more dominant SVD components.

Using the definition of canonical angles $\Theta(\mathcal{X},\mathcal{Y})$
between the two subspaces $\mathcal{X}$ and $\mathcal{Y}$ of the same dimension
\cite[p. 250]{stewart01}, we have the following theorem, which
shows how well the subspace $\mathcal{K}_{k}(A^{T}A,A^{T}b)$, on which
LSQR and CGLS work, captures the $k$-dimensional dominant right singular space.

\begin{theorem}\label{thm2}
Let the SVD of $A$ be \eqref{eqsvd}, and assume that its singular values
are distinct and satisfy $\sigma_j=\mathcal{O}(e^{-\alpha j})$ with $\alpha>0$.
Let $\mathcal{V}_k=span\{V_k\}$ be the subspace spanned by the columns of
$V_k=(v_1,v_2,\ldots,v_k)$, and $\mathcal{V}_k^s=\mathcal{K}_{k}
(A^{T}A, A^{T}b)$. Then
\begin{equation}\label{sindeltak}
\|\sin\Theta(\mathcal{V}_k,\mathcal{V}_k^s)\|=\frac{\|\Delta_k\|}
{\sqrt{1+\|\Delta_k\|^2}}
\end{equation}
with the $(n-k)\times k$ matrix $\Delta_k$ to be defined by
\eqref{deltakdef} and
\begin{equation}\label{eqres1}
  \|\Delta_k\|_F\leq
  \frac{\sigma_{k+1}}{\sigma_k}\frac{\max_{j=k+1}^n|u_{j}^Tb|}{\min_{j=1}^k|u_{j}^Tb|}
  \sqrt{k(n-k)}(1+\mathcal{O}(e^{-2\alpha})), \ k=1,2,\ldots,n-1.
\end{equation}
\end{theorem}

{\em Proof}.
Let $\bar{U}=(u_1,u_2,\ldots,u_n)$ consist of the first $n$ columns of
$U$ defined in \eqref{eqsvd}. We see $\mathcal{K}_{k}(\Sigma^2,
\Sigma \bar{U}^Tb)$ is spanned by the columns of the $n\times k$ matrix $DT_k$
with
\begin{equation*}\label{}
  D=\mathrm{diag}\left(\sigma_i \bar{U}^Tb\right),\ \ \
  T_k=\left(\begin{array}{cccc} 1 &
  \sigma_1^2&\ldots & \sigma_1^{2k-2}\\
1 &\sigma_2^2 &\ldots &\sigma_2^{2k-2} \\
\vdots & \vdots&&\vdots\\
1 &\sigma_n^2 &\ldots &\sigma_n^{2k-2}
\end{array}\right).
\end{equation*}
Partition the matrices $D$ and $T_k$ as follows:
\begin{equation*}\label{}
  D=\left(\begin{array}{cc} D_1 & 0 \\ 0 & D_2 \end{array}\right),\ \ \
  T_k=\left(\begin{array}{c} T_{k1} \\ T_{k2} \end{array}\right),
\end{equation*}
where $D_1, T_{k1}\in\mathbb{R}^{k\times k}$. Since $T_{k1}$ is
a Vandermonde matrix with $\sigma_j$ distinct for $1\leq j\leq k$, it is
nonsingular. Thus, by the SVD of $A$, we have
\begin{equation*}\label{}
  \mathcal{K}_{k}(A^{T}A, A^{T}b)=span\{VDT_k\}=span
  \left\{V\left(\begin{array}{c} D_1T_{k1} \\ D_2T_{k2} \end{array}\right)\right\}
  =span\left\{V\left(\begin{array}{c} I\\ \Delta_k \end{array}\right)\right\}
\end{equation*}
with
\begin{equation}\label{deltakdef}
\Delta_k=D_2T_{k2}T_{k1}^{-1}D_1^{-1}.
\end{equation}
Define $Z_k=V\left(\begin{array}{c} I \\ \Delta_k \end{array}\right)$.
Then $Z_k^TZ_k=I+\Delta_k^T\Delta_k$ and the columns of
$Z_k(Z_k^TZ_k)^{-\frac{1}{2}}$ form an orthonormal basis of $\mathcal{V}_k^s$.

Write $V=(V_k, V_k^{\perp})$. By definition, we obtain
\begin{align}
   \|\sin\Theta(\mathcal{V}_k,\mathcal{V}_k^s)\|  &=
   \left\|(V_k^{\perp})^T Z_k(Z_k^TZ_k)^{-\frac{1}{2}}\right\|\notag\\
   &=\left\|(V_k^{\perp})^T V
   \left(\begin{array}{c} I \\ \Delta_k \end{array}\right)
   (I+\Delta_k^T\Delta_k)^{-\frac{1}{2}}\right\| \notag\\
   &=\|\Delta_k (I+\Delta_k^T\Delta_k)^{-\frac{1}{2}}\|
   =\frac{\|\Delta_k\|}{\sqrt{1+\|\Delta_k\|^2}},\notag
\end{align}
which proves \eqref{sindeltak} and indicates that
$\|\sin\Theta(\mathcal{V}_k,\mathcal{V}_k^s)\|$ is monotonically increasing
with respect to $\|\Delta_k\|$.

We next estimate $\|\Delta_k\|$.  We have
\begin{align}
\|\Delta_k\|&\leq\|\Delta\|_F=\left\|D_2T_{k2}T_{k1}^{-1}D_1^{-1}\right\|_F
   \leq \|D_2\|\left\|T_{k2}T_{k1}^{-1}\right\|_F\left\|D_1^{-1}\right\|\notag\\
 &=\frac{\sigma_{k+1}}{\sigma_k}\frac{\max_{j=k+1}^n|v_{j}^Tb|}
 {\min_{j=1}^k|v_{j}^Tb|}\left\|T_{k2}T_{k1}^{-1}\right\|_F.
   \label{fnorm}
\end{align}
So we need to estimate $\left\|T_{k2}T_{k1}^{-1}\right\|_F$.
It is easily justified that the $i$-th column of $T_{k1}^{-1}$ consists of
the coefficients of the Lagrange polynomial
\begin{equation*}\label{}
  L_i^{(k)}(\lambda)=\prod\limits_{j=1,j\neq i}^k
  \frac{\sigma_j^2-\lambda}{\sigma_j^2-\sigma_i^2}
\end{equation*}
that interpolates the elements of the $i$-th canonical basis vector
$e_i^{(k)}\in \mathbb{R}^{k}$ at the abscissas $\sigma_1^2,\sigma_2,
\ldots, \sigma_k^2$. Consequently, the $i$-th column of $T_{k2}T_{k1}^{-1}$ is
\begin{equation*}\label{}
  T_{k2}T_{k1}^{-1}e_i^{(k)}=\left(L_i^{(k)}(\sigma_{k+1}^2),
  \ldots,L_i^{(k)}(\sigma_{n}^2)\right)^T,
\end{equation*}
from which we obtain
\begin{equation}\label{tk2tk1}
  T_{k2}T_{k1}^{-1}=\left(\begin{array}{cccc} L_1^{(k)}(\sigma_{k+1}^2)&
  L_2^{(k)}(\sigma_{k+1}^2)&\ldots & L_k^{(k)}(\sigma_{k+1}^2)\\
L_1^{(k)}(\sigma_{k+2}^2)&L_2^{(k)}(\sigma_{k+2}^2) &\ldots &
L_k^{(k)}(\sigma_{k+2}^2) \\
\vdots & \vdots&&\vdots\\
L_1^{(k)}(\sigma_{n}^2)&L_2^{(k)}(\sigma_{n}^2) &\ldots &L_k^{(k)}(\sigma_{n}^2)
\end{array}\right).
\end{equation}
Since $|L_i^{(k)}(\lambda)|$ is monotonic for $\lambda<\sigma_k^2$,
it is bounded by $|L_i^{(k)}(0)|$.
Furthermore, let $|L_{i_0}^{(k)}(0)|=\max_{i=1,2,\ldots,k}|L_i^{(k)}(0)|$.
Then for $i=1,2,\ldots,k$ and $\alpha>0$ we have
\begin{align}
|L_i^{(k)}(0)| &\leq |L_{i_0}^{(k)}(0)|=\prod\limits_{j=1,j\neq i_0}^k
\left|\frac{\sigma_j^2}{\sigma_j^2-\sigma_{i_0}^2}\right|
  =\prod\limits_{j=1}^{i_0-1}\frac{\sigma_j^2}{\sigma_j^2-\sigma_{i_0}^2}
   \cdot\prod\limits_{j=i_0+1}^{k}\frac{\sigma_j^2}{\sigma_{i_0}^2-\sigma_{j}^2} \notag\\
& =\prod\limits_{j=1}^{i_0-1}\frac{1}
{1-\mathcal{O}(e^{-2(i_0-j)\alpha})}
\prod\limits_{j=i_0+1}^{k}\frac{1}
{\mathcal{O}(e^{2(j-i_0)\alpha})-1} \notag\\
& =\prod\limits_{j=1}^{i_0-1}\frac{1}
{1-\mathcal{O}(e^{-2(i_0-j)\alpha})}
\prod\limits_{j=i_0+1}^{k}\frac{1}
{1-\mathcal{O}(e^{-2(j-i_0)\alpha})}\frac{1}
{\prod\limits_{j=i_0+1}^{k}\mathcal{O}(e^{2(j-i_0)\alpha})} \notag\\
&=\frac{\left(1+\sum\limits_{j=1}^{i_0} \mathcal{O}(e^{-2j\alpha})\right)
\left(1+\sum\limits_{j=1}^{k-i_0+1} \mathcal{O}(e^{-2j\alpha})\right)}
{\prod\limits_{j=i_0+1}^{k}\mathcal{O}(e^{2(j-i_0)\alpha})} \label{lik}
\end{align}
by absorbing those higher order terms into ${\cal O}(\cdot)$.
 Note that in the above numerator we have
  $$
  1+\sum\limits_{j=1}^{i_0} \mathcal{O}(e^{-2j\alpha})
    =1+ \mathcal{O}\left(\sum\limits_{j=1}^{i_0}e^{-2j\alpha}\right)
    =1+ \mathcal{O}\left(\frac{e^{-2\alpha}}
    {1-e^{-2\alpha}}(1-e^{-2i_0\alpha})\right),
  $$
  and
  $$
    1+\sum\limits_{j=1}^{k-i_0+1} \mathcal{O}(e^{-2j\alpha})
    =1+ \mathcal{O}\left(\sum\limits_{j=1}^{k-i_0+1}e^{-2j\alpha}\right)
    =1+ \mathcal{O}\left(\frac{e^{-2\alpha}}{1-e^{-2\alpha}}
    (1-e^{-2(k-i_0+1)\alpha})\right).
  $$
It is then easily seen that their product is
  $$
  1+ \mathcal{O}\left(\frac{2e^{-2\alpha}}{1-e^{-2\alpha}}\right)
  +\mathcal{O}\left(\left(\frac{e^{-2\alpha}}{1-e^{-2\alpha}}\right)^2\right)=
  1+ \mathcal{O}\left(\frac{2e^{-2\alpha}}{1-e^{-2\alpha}}\right)
  =  1+\mathcal{O}(e^{-2\alpha}).
  $$
On the other hand, by definition, the denominator
  $\prod\limits_{j=i_0+1}^k\left(\frac{\sigma_{i_0}}{\sigma_j}\right)^2
=\prod\limits_{j=i_0+1}^{k}\mathcal{O}(e^{2(j-i_0)\alpha})$ in
  \eqref{lik} is exactly one for $i_0=k$, and it
  is strictly bigger than one for $i_0<k$. Therefore, for any $k$,
  we have
  $|L_{i_0}^{(k)}(0)|=\max_{i=1,2,\ldots,k} |L_i^{(k)}(0)|=1+
  \mathcal{O}(e^{-2\alpha})$.
From this and \eqref{tk2tk1} it follows
that
\begin{equation*}\label{}
\left\|T_{k2}T_{k1}^{-1}\right\|_F\leq \sqrt{k}
\left\|T_{k2}T_{k1}^{-1}e_k^{(k)}\right\|\leq
\sqrt{k (n-k)}|L_{i_0}^{(k)}(0)|
=\sqrt{k (n-k)}(1+\mathcal{O}(e^{-2\alpha})).
\end{equation*}
Therefore, for $i=1,2,\ldots,n-1$  and $\alpha>0$ considerably, from
\eqref{fnorm} we have
\begin{align}\label{b1b2}
 \|\Delta_k\|_F &\leq\frac{\sigma_{k+1}}{\sigma_{k}}
 \frac{\max_{j=k+1}^n|u_{j}^Tb|}{\min_{j=1}^k|u_{j}^Tb|}
 \sqrt{k(n-k)}(1+\mathcal{O}(e^{-2\alpha})). \qquad\endproof
\end{align}

{\bf Remark 2.1}\ \ \ We point out that \eqref{eqres1} should not be sharp.
As we have seen from the proof, the factor $\frac{\sigma_{k+1}}
{\sigma_k}\frac{\max_{j=k+1}^n|u_{j}^Tb|}{\min_{j=1}^k|u_{j}^Tb|}$
seems intrinsic and unavoidable,
but the factor $\sqrt{k(n-k)}$ in \eqref{eqres1} is an overestimate
and can certainly be reduced. \eqref{b1b2} is an overestimate
since $|L_i^{(k)}(0)|$ for $i$ not near to $k$ is considerably
smaller than $|L_{i_0}^{(k)}(0|$, but we replace all them by their
maximum $1+\mathcal{O}(e^{-2\alpha})$. In fact, our derivation
clearly illustrates that the smaller $i$ is, the smaller
$|L_i^{(k)}(0)|$ than $|L_k^{(k)}(0)|$.

Recall the discrete Picard condition \eqref{picard}. Then
\begin{equation}\label{ck}
  c_k=\frac{\max_{j=k+1}^n|u_{j}^Tb|}{\min_{j=1}^k|u_{j}^Tb|}=
  \frac{\max_{j=k+1}^n(|u_{j}^T\hat{b}+u_{j}^Te|)}
  {\min_{j=1}^k(|u_{j}^T\hat{b}+u_{j}^Te|)}\approx
\frac{\sigma_{k+1}^{1+\beta}+|u_{k+1}^Te|}{\sigma_{k}^{1+\beta}+|u_{k}^Te|}.
\end{equation}
We observe that $c_k\approx \frac{\sigma_{k+1}^{1+\beta}}{\sigma_k^{1+\beta}}<1$
almost remains constant for $k\leq k_0$.
For $k>k_0$, note that all the $|u_k^T b|\approx |u_{k}^Te|$ almost remain the
same. Thus, we have $c_k\approx 1$, meaning that $\mathcal{V}_k^s$ does not
capture $\mathcal{V}_k$ as well as it does for $k\leq k_0$.

{\bf Remark 2.2}\ \ \ The theorem can be extended to
moderately ill-posed problems with the singular values
$\sigma_j=\mathcal{O}(j^{-\alpha}),\ \alpha>1$ considerably and $k$
not big since, in a similar manner to the proof of Theorem~\ref{thm2},
we can obtain by the first order Taylor expansion
\begin{align*}\label{}
  |L_{i_0}^{(k)}(0)|&\approx |L_k^{(k)}(0)|= \prod\limits_{j=1}^{k-1}
\frac{\sigma_j^2}{\sigma_j^2-\sigma_k^2}\\
&=\prod\limits_{j=1}^{k-1}\frac{1}
{1-\mathcal{O}((\frac{j}{k})^{2\alpha})}\approx 1+
\sum\limits_{j=1}^{k-1}\mathcal{O}\left(\left(\frac{j}{k}\right)^
{2\alpha}\right)
=\mathcal{O}(1),
\end{align*}
which, unlike $(1+\mathcal{O}(e^{-2\alpha}))$ for severely ill-posed problems,
depends on $k$ and increases slowly with $k$ for $\alpha>1$ considerably.
However, for mildly ill-posed problems,
from above we have $|L_{i_0}^{(k)}(0)|>1$ considerably for $\alpha<1$.

{\bf Remark 2.3}\ \ A combination of \eqref{sindeltak} and \eqref{eqres1}
and the above analysis indicate that $\mathcal{V}_k^s$  captures
$\mathcal{V}_k$ better for severely ill-posed problems than for moderately
ill-posed problems. There are two reasons for this.
The first is that the factors $\sigma_{k+1}/\sigma_k$ are
basically fixed constants for severely ill-posed problems as $k$
increases, and they are smaller than the counterparts for moderately ill-posed
problems unless the degree $\alpha$ of its ill-posedness is far bigger than one
and $k$ small. The second is that the factor $\mathcal{O}(1)$
is smaller for severely ill-posed problems than the factor
$1+\mathcal{O}(e^{-2\alpha})$ for moderately ill-posed problems for the same $k$.

{\bf Remark 2.4} \ \
The situation is fundamentally different for mildly ill-posed problems:
Firstly, we always have $|L_{i_0}^{(k)}(0)|>1$ substantially for
$\alpha\leq 1$ and any $k$, which is considerably bigger than $\mathcal{O}(1)$
for moderately ill-posed problems for the same $k$. Secondly,
$c_k$ defined by \eqref{ck} is closer to one than that
for moderately ill-posed problems for $k=1,2,\ldots,k_0$.
Thirdly, for the same noise level $\|e\|$ and $\beta$, we see from the discrete
Picard condition \eqref{picard} and the definition of $k_0$ that $k_0$
is bigger for a mildly ill-posed problem than that for a moderately
ill-posed problem. All of them show that $\mathcal{V}_k^s$
captures $\mathcal{V}_k$ {\em considerably better} for severely and moderately
ill-posed problems than for mildly ill-posed problems for $k=1,2,\ldots,k_0$.
In other words, our results illustrate that $\mathcal{V}_k^s$ contains more
information on the other $n-k$ right singular vectors for mildly ill-posed
problems, compared with severely and moderately ill-posed problems.
The bigger $k$, the more it contains. Therefore, $\mathcal{V}_k^s$ captures
$\mathcal{V}_k$ more effectively for severely and moderately ill-posed
problems than  mildly ill-posed problems. That is, $\mathcal{V}_k^s$
contains more information on the other $n-k$ right singular vectors for mildly
ill-posed problems, making the appearance of a small Ritz value more possible
before $k\leq k_0$ and LSQR has better regularizing effects
for the first two kinds of problems than for the third kind.
Note that LSQR, at most, has the full regularization, i.e., there is no Ritz
value smaller than $\sigma_{k_0+1}$ for $k\leq k_0$, for severely and
moderately ill-posed problems. Our analysis indicates that LSQR generally
has only the partial regularization for mildly ill-posed problem and a hybrid
LSQR should be used.

{\bf Remark 2.5}\ \ \ Relation \eqref{eqres1} and $c_k$ indicate
that $\mathcal{V}_k^s$ captures $\mathcal{V}_k$ better for severely ill-posed
problems than for moderately ill-posed problems. There are two reasons for this.
First, the all the $\sigma_{k+1}/\sigma_k$ are basically a fixed constant
$\rho^{-1}$ for severely ill-posed problems, which is smaller than those
ratios for moderately
ill-posed problems unless $\alpha$ is rather big and $k$ small. Second,
the quantities $|L_{i_0}^{(k)}(0)|=1+\mathcal{O}(e^{-2\alpha})$
for severely ill-posed problems are smaller than the corresponding
$\mathcal{O}(1)$ for moderately ill-posed problems.

Let us investigate more and get insight into the regularization of LSQR. Define
\begin{equation}\label{gammak}
\gamma_k = \left\|A-P_{k+1}B_kQ_k^T\right\|,
\end{equation}
which measures the quality of the rank $k$ approximation $P_{k+1}B_kQ_k^T$ to $A$.
Based on \eqref{eqres1}, we can derive the following estimate for $\gamma_k$.

\begin{theorem}\label{thm3}
Assume that \eqref{eq1} is severely or moderately ill posed. Then
\begin{equation}\label{gamma}
\sigma_{k+1}\leq  \gamma_k \leq  \sigma_{k+1} +
\sigma_1\|\sin\Theta(\mathcal{V}_k,\mathcal{V}_k^s)\|.
\end{equation}
\end{theorem}

{\em Proof}.
Let $A_k=U_k\Sigma_k V_k^T$ be the best rank $k$ approximation to $A$
with respect to the 2-norm, where $U_k=(u_1,\ldots,u_k)$,
$V_k=(v_1,\ldots,v_k)$ and $\Sigma_k={\rm diag}(\sigma_1,\ldots,\sigma_k)$.
Since $P_{k+1}B_kQ_k^T$ is of rank $k$, the lower bound in \eqref{gamma}
is trivial by noting that $\gamma_k\geq \|A-A_k\|=\sigma_{k+1}$.
We now prove the upper bound. From \eqref{eqmform1}, we obtain
\begin{align*}
\left\|A-P_{k+1}B_kQ_k^T\right\|&= \left\|A-AQ_kQ_k^T\right\|.
\end{align*}
It is easily known that $\mathcal{V}_k^s
=\mathcal{K}_{k}(A^{T}A,A^{T}b)=span\{Q_k\}$ with $Q_k$ having orthonormal
columns. Then by the definition of $\|\sin\Theta(\mathcal{V}_k,\mathcal{V}_k^s)\|$
we obtain
\begin{align*}
 \left\|A-AQ_kQ_k^T\right\| &= \left\|(A-U_k\Sigma_kV_k^T+U_k\Sigma_kV_k^T)
     (I-Q_kQ_k^T)\right\| \notag\\
 &\leq\left\|(A-U_k\Sigma_kV_k^T)(I-Q_kQ_k^T)\right\|+
     \left\|U_k\Sigma_kV_k^T(I-Q_kQ_k^T)\right\| \notag\\
 &\leq \sigma_{k+1} + \|\Sigma_k\|\left\|V_k^T(I-Q_kQ_k^T)\right\| \notag\\
 &= \sigma_{k+1} + \sigma_1\|\sin\Theta(\mathcal{V}_k,\mathcal{V}_k^s)\|.
 \qquad\endproof
\end{align*}

Numerically, it has been extensively observed in the literature that the
$\gamma_k$ decay as fast as $\sigma_{k+1}$ and, more precisely,
$\gamma_k\approx \sigma_{k+1}$ for severely ill-posed problems;
see, e.g., \cite{bazan14,gazzola14,gazzola-online}. They mean that
the $P_{k+1}B_kQ_k^T$ are very good rank $k$ approximations to $A$.
Recall that the TSVD method generates the best regularized solution
$x_{k_0}^{TSVD}=A_{k_0}^{\dagger}b$. As a result, if $\gamma_{k_0}\approx
\sigma_{k_0+1}$, the LSQR iterate $x^{(k_0)}=Q_{k_0}T_{k_0}^{\dagger}
Q_{k_0+1}^Tb$ is reasonably close to the TSVD solution $x_{k_0}^{TSVD}$
for $\sigma_{k_0+1}$ is reasonably small. This means that LSQR
has the full regularization and does not need any additional regularization
to improve $x^{(k_0)}$.
As our experiments will indicate in detail, these observed phenomena are of
generality for both severely and moderately ill-posed problems
and thus should have strong theoretical supports. Compared to the observations,
our \eqref{gamma} appears to be a considerrable overestimate.

We next present some results on $\alpha_{k+1}$ appearing in \eqref{eqmform2}.
If $\alpha_{k+1}=0$, the Lanczos bidiagonalization process terminates,
and we have found $k$ exact singular triples of $A$ \cite{jia03}. In our
context, since $A$ has only simple singular values and $b$ has components in all
the left singular vectors, early termination is impossible in exact arithmetic,
but small $\alpha_{k+1}$ is possible.
We aim to investigate how fast $\alpha_{k+1}$ decays.
We first give a refinement of a result in \cite{gazzola-online}.

\begin{theorem}\label{thm1}
Let $B_k=W_k{\Theta_k}S_k^{T}$ be the SVD of $B_k$,
where $W_k\in\mathbb{R}^{(k+1)\times (k+1)}$ and
$S_k\in\mathbb{R}^{k\times k}$ are orthogonal, and
$\Theta_k
 \in\mathbb{R}^{(k+1)\times k}$,
and define $\tilde{U}_k=P_{k+1}{W_k}$ and $\tilde{V}_k=Q_{k}{S_k}$. Then
\begin{eqnarray}
  A\tilde{V}_k-\tilde{U}_k{\Theta_k}&=&0,\label{refine}\\
  \left\|A^{T}\tilde{U}_k-\tilde{V}_k\Theta_k^T\right\|&=&
  \alpha_{k+1}.\label{eqpro1}
\end{eqnarray}
\end{theorem}

{\em Proof}.
From \eqref{eqmform1} and $B_k=W_k{\Theta_k}S_k^T$, we obtain
\begin{equation*}\label{}
  A\tilde{V}_k=AQ_kS_k=P_{k+1}B_kS_k=\tilde{U}_k{\Theta_k}.
\end{equation*}
So \eqref{refine} holds. From \eqref{eqmform2}, we get
\begin{align*}
A^{T}\tilde{U}_k  &=  A^{T}P_{k+1}{W_k} \notag\\
&= Q_kB_k^{T}W_k+ \alpha_{k+1}q_{k+1}e_{k+1}^{T}W_k\notag\\
&= Q_kS_k\Theta_k^{T}+ \alpha_{k+1}q_{k+1}e_{k+1}^{T}W_k \notag\\
&= \tilde{V}_k\Theta_k^{T}+ \alpha_{k+1}q_{k+1}e_{k+1}^{T}W_k.
\end{align*}
Note that $\|q_{k+1}\|=\|e_{k+1}^TW_k\|=1$. Then we get
\begin{equation}\label{eqalpha}
  \left\|A^{T}\tilde{U}_k-\tilde{V}_k\Theta_k^T\right\|=\alpha_{k+1}
  \left\|q_{k+1}e_{k+1}^TW_k\right\|=\alpha_{k+1}. \qquad\endproof
\end{equation}

We remark that it is an inequality other than the equality
in a result of \cite{gazzola-online} similar to \eqref{eqpro1}.

In combination with the previous results and remarks,
this theorem shows that once $\alpha_{k+1}$ becomes
small for not big $k$, the $k$ singular values of $B_k$ may
approximate the large singular values of $A$, and it is more possible
that no small one appears for severely ill-posed problems and
moderately ill-posed problems.

As our final result, we establish an intimate and interesting relationship
between $\alpha_{k+1}$ and $\gamma_k$, showing how fast $\alpha_{k+1}$ decays.

\begin{theorem}\label{thm4}
It holds that
\begin{equation}\label{}
  \alpha_{k+1}\leq\gamma_k.
\end{equation}
\end{theorem}

{\em Proof}.
With the notations as in Theorem~\ref{thm1},
we have $P_{k+1}B_kQ_k^T=\tilde{U}_k\Theta_k\tilde{V}_k^T$.
So, by \eqref{gammak}, we have
$$
\gamma_k=\left\|A-\tilde{U}_k\Theta_k\tilde{V}_k^T\right\|.
$$
Note that $\tilde{U}_k^T\tilde{U}_k=I$. Therefore, from \eqref{eqalpha} we obtain
\begin{align*}
\alpha_{k+1} &=\left\|A^{T}\tilde{U}_k-\tilde{V}_k\Theta_k^T\right\| \notag\\
           &=\left\|A^{T}\tilde{U}_k\tilde{U}_k^T-
           \tilde{V}_k\Theta_k^T\tilde{U}_k^T\right\| \notag\\
           &=\left\|A^{T}\tilde{U}_k\tilde{U}_k^T-\tilde{V}_k\Theta_k^T\tilde{U}_k^T
           \tilde{U}_k\tilde{U}_k^T\right\|\notag\\
           &=\left\|\left(A^{T}-\tilde{V}_k\Theta_k^T\tilde{U}_k^T\right)
           \tilde{U}_k\tilde{U}_k^T\right\| \notag\\
           &\leq \left\|\left(A^{T}-\tilde{V}_k\Theta_k^T\tilde{U}_k^T\right)\right\|
           \left\|\tilde{U}_k\tilde{U}_k^T\right\| \notag\\
           &=\left\|A-\tilde{U}_k\Theta_k\tilde{V}_k^T\right\| = \gamma_k. \qquad\endproof
\end{align*}

The theorem indicates that $\alpha_{k+1}$ decays at least as fast as
$\gamma_k$, which, in turn, means that $\alpha_{k+1}$ may decrease
in the same rate as $\sigma_{k+1}$, as observed
in \cite{bazan14,gazzola14,gazzola-online} for severely
ill-posed problems.

\section{Numerical experiments}\label{SectionExp}

In this section, we report numerical experiments to illustrate the
the regularizing effects of LSQR. We will demonstrate that LSQR has
the full regularization for severely and moderately ill-posed problems,
stronger phenomena than our theory proves, but it only
has the partial regularization for mildly ill-posed problems, in
accordance with our theory, for which a hybrid LSQR is needed to
compute best possible regularized solutions. We choose several
ill-posed examples from Hansen's regularization toolbox \cite{hansen07}.
All the problems arise from the discretization of the first kind Fredholm
integral equation
\begin{equation}\label{eq2}
  \int_a^bK(s,t)x(t)dt=b(s),\ \ \ c\leq s \leq d.
\end{equation}
For each problem we use the codes of \cite{hansen07} to
generate a $1024\times 1024$ matrix $A$, true solution $x_{true}$
and noise-free right-hand $\hat{b}$. In order to simulate the noisy data,
we generate the Gaussian noise vector $e$ whose entries are normally
distributed with mean zero. Defining the noise level
$\varepsilon=\frac{\|e\|}{\|\hat{b}\|}$, we use $\varepsilon=10^{-2},
10^{-3}, 10^{-4}$, respectively, in the test examples.
To simulate exact arithmetic, the full reorthogonalization is used
during the Lanczos bidiagonalization process.
We remind that, as far as ill-posed problem \eqref{eq1}
is concerned,
our primary goal consists in justifying the regularizing effects
of iterative solvers, which are {\em unaffected by sizes}
of ill-posed problems and only depends on the degree of
ill-posedness. Therefore, for this purpose, as extensively done in the
literature (see, e.g., \cite{hansen98,hansen10} and the references therein),
it is enough to test not very large problems. Indeed, for
$n$ large, say, 1,0000 and more, we have observed completely the same
behavior as that for $n$ not large, e.g., $n=1024$ used in this paper.
A reason for using $n$ not large is because
such choice makes it practical to fully
justify the regularization effects of LSQR by comparing it with
the TSVD method, which suits only for
small and/or medium sized problems for computational efficiency.
All the computations are carried out in Matlab 7.8 with the
machine precision $\epsilon_{\rm mach}=2.22\times10^{-16}$ under the
Microsoft Windows 7 64-bit system.

\subsection{Severely ill-posed problems}
We consider the following two severely ill-posed problems \cite{hansen07}.

{\bf Example 1}\ \ \ This problem 'Shaw' arises from one-dimensional
image restoration, and can be obtained by discretizing the first kind
Fredholm integral equation \eqref{eq2} with $[-\frac{\pi}{2}, \frac{\pi}{2}]$ as
both integration and domain intervals. The kernel $K(s,t)$ and the solution $x(t)$
are given by
\begin{equation*}\label{}
  K(s,t)=(\cos(s)+\cos(t))^2\left(\frac{\sin(u)}{u}\right)^2,\ \ \
  u=\pi(\sin(s)+\sin(t)),
\end{equation*}
\begin{equation*}\label{}
  x(t)=2\exp(-6(t-0.8)^2)+\exp(-2(t+0.5)^2).
\end{equation*}

{\bf Example 2}\ \ \ This problem 'Wing' has a discontinuous solution
and is obtained by discretizing the first kind Fredholm integral
equation \eqref{eq2} with $[0, 1]$ as both integration and domain intervals.
The kernel $K(s,t)$, the solution $x(t)$ and the right-hand side $b(s)$
are given by
\begin{equation*}\label{}
  K(s,t)=t\exp(-st^2),\ \ \ b(s)=\frac{\exp(-\frac{1}{9}s)-
  \exp(-\frac{4}{9}s)}{2s},
\end{equation*}
\begin{equation*}\label{}
  x(t)=\left\{\begin{array}{ll} 1,\ \ \ &\frac{1}{3}<t<\frac{2}{3};\\ 0,\ \ \
  &elsewhere. \end{array}\right.
\end{equation*}

These two problems are severely ill-posed, whose singular values
$\sigma_j=\mathcal{O}(e^{-\alpha j})$ with $\alpha=2$ for 'Shaw' and
$\alpha=4.5$ for 'Wing', respectively.

\begin{figure}

\begin{minipage}{0.48\linewidth}
  \centerline{\includegraphics[width=7.0cm,height=5cm]{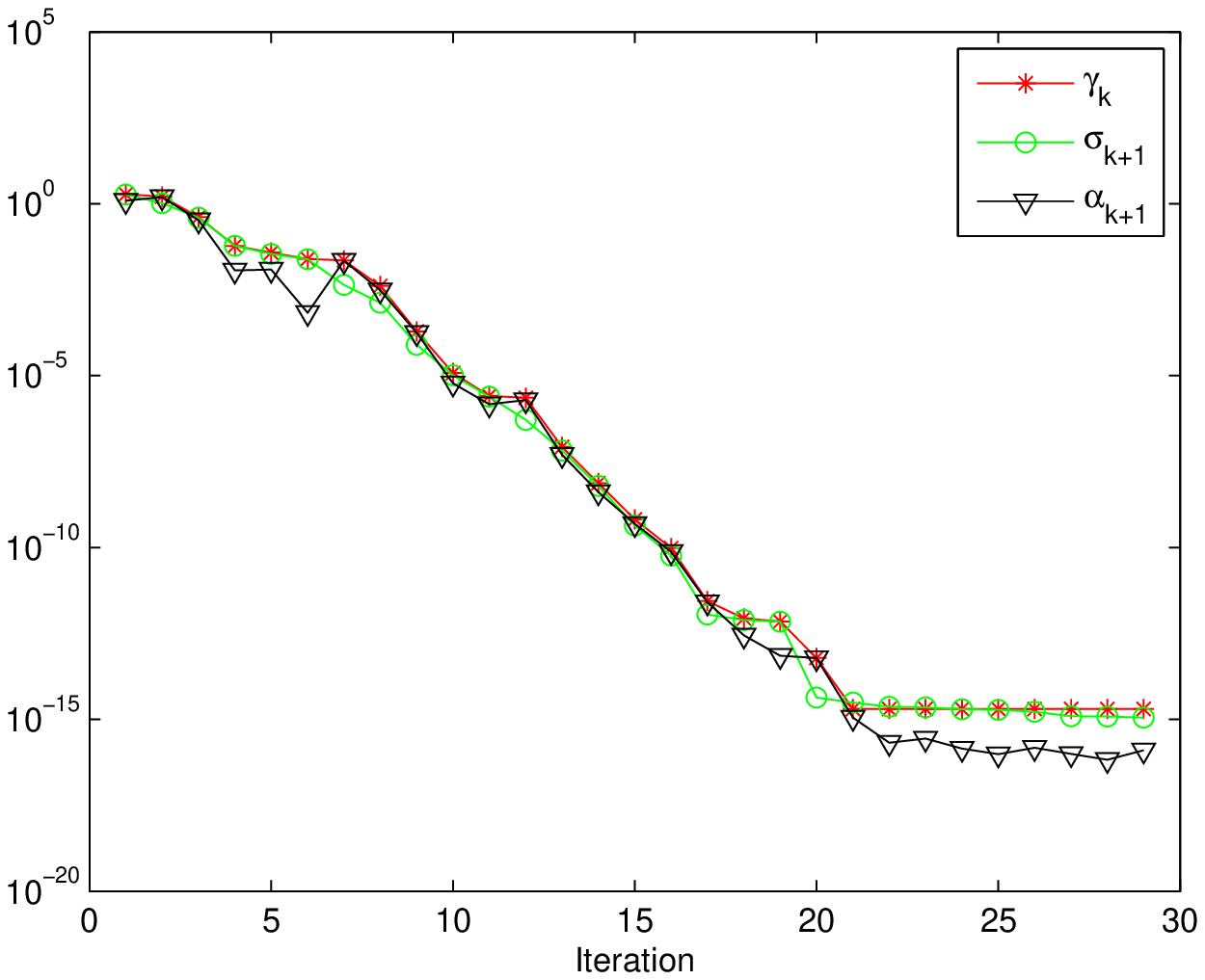}}
  \centerline{(a)}
\end{minipage}
\hfill
\begin{minipage}{0.48\linewidth}
  \centerline{\includegraphics[width=7.0cm,height=5cm]{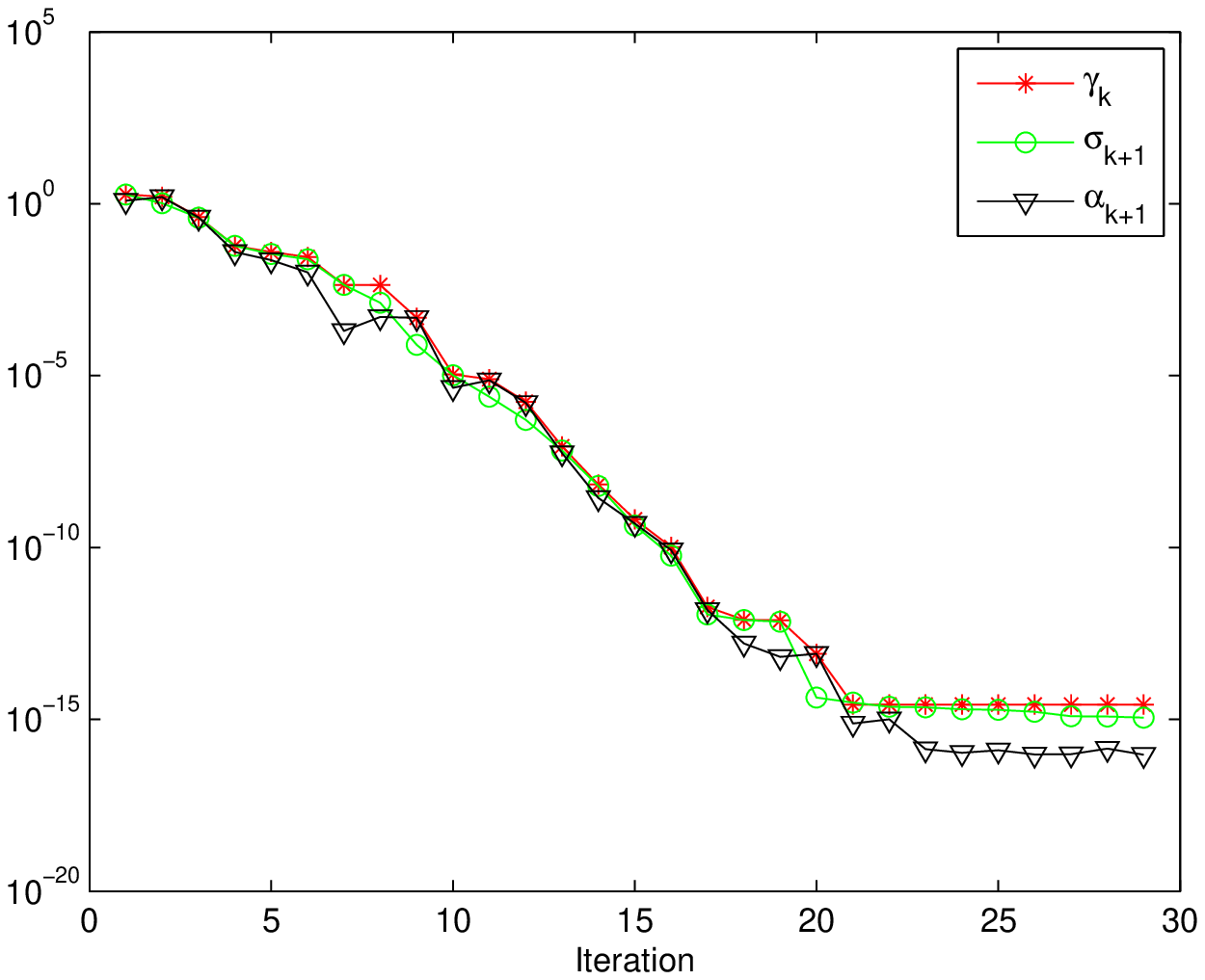}}
  \centerline{(b)}
\end{minipage}
\vfill
\begin{minipage}{0.48\linewidth}
  \centerline{\includegraphics[width=7.0cm,height=5cm]{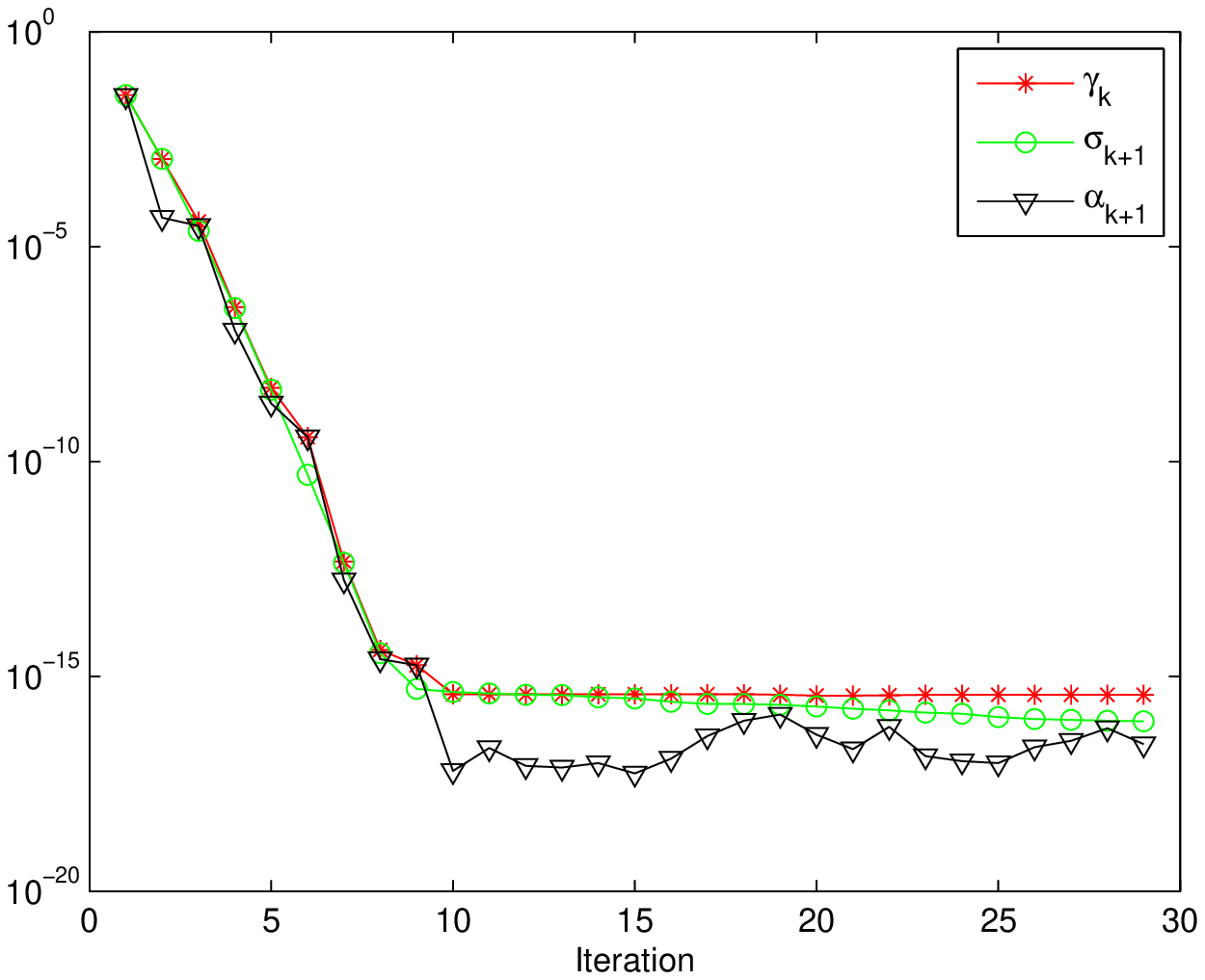}}
  \centerline{(c)}
\end{minipage}
\hfill
\begin{minipage}{0.48\linewidth}
  \centerline{\includegraphics[width=7.0cm,height=5cm]{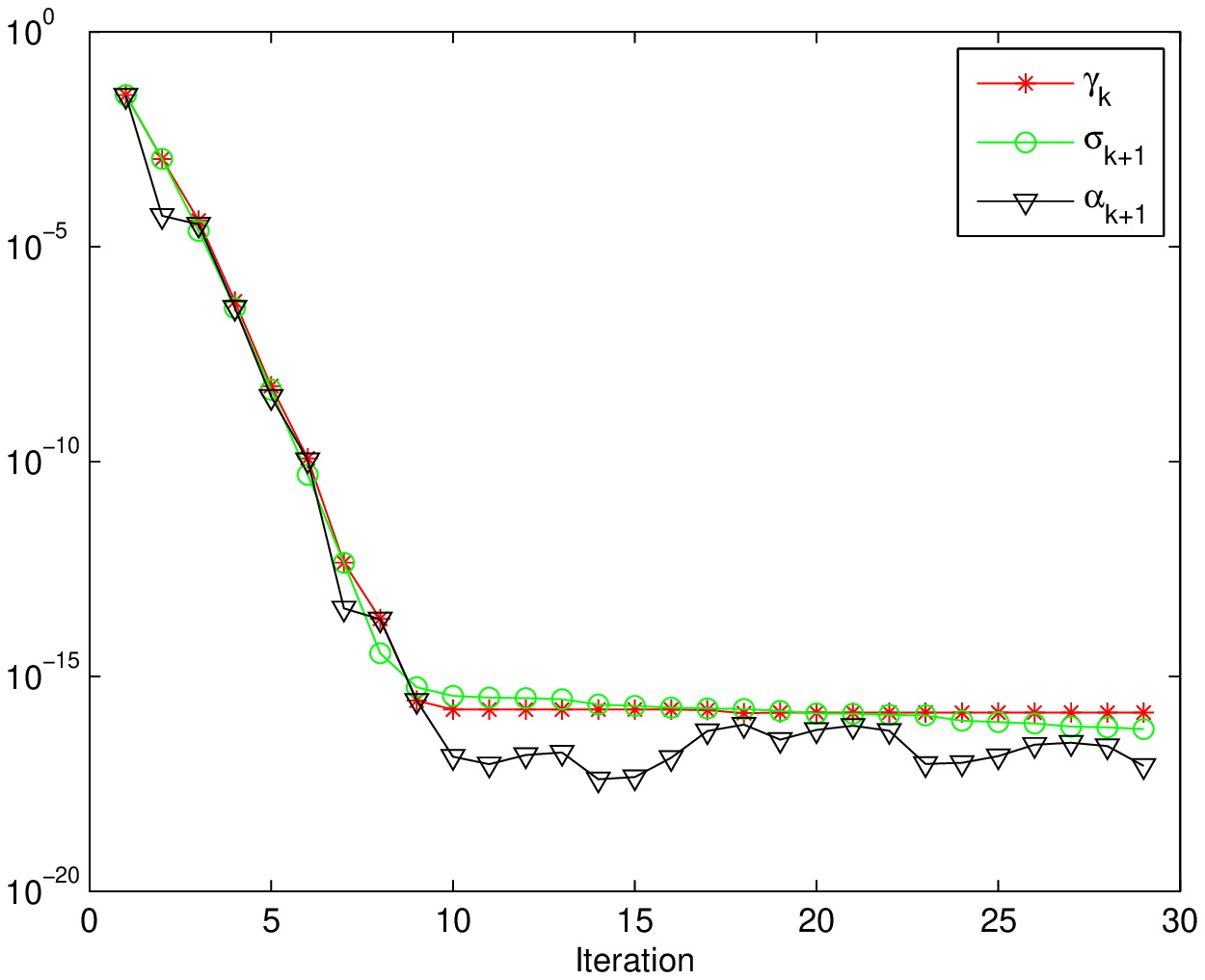}}
  \centerline{(d)}
\end{minipage}
\caption{(a)-(b): Plots of decaying behavior of the sequences $\gamma_k$, $\sigma_{k+1}$ and
$\alpha_{k+1}$ for the problem Shaw with $\varepsilon=10^{-2}$ (left) and
$\varepsilon=10^{-3}$ (right); (c)-(d): Plots of decaying behavior of the
sequences $\gamma_k$ and $\sigma_{k+1}$ for the problem Wing with
$\varepsilon=10^{-3}$ (left) and $\varepsilon=10^{-4}$ (right).}
\label{fig:res}
\end{figure}

In Figure~\ref{fig:res}, we display the curves of the sequences $\gamma_k$ and
$\alpha_{k+1}$ with $\varepsilon=10^{-2}, 10^{-3}, 10^{-4}$, respectively.
They illustrate that the quantities $\gamma_k$ decrease as fast as $\sigma_{k+1}$
and both of them level off at the level of $\epsilon_{\rm mach}$ for $k$ no more
than 20, and after that these quantities are purely round-offs and
are reliable no more. Moreover, the curves of quantities $\alpha_{k+1}$ always
lie below those of $\gamma_k$, which coincides with Theorem \ref{thm4}.
We can see that the decaying curves with different noise levels are almost
the same. Furthermore, we observe that
$\gamma_k\approx \sigma_{k+1}$ for severely ill-posed problems,
indicating that the $P_{k+1}B_kQ_k^T$ are very good rank $k$ approximations
to $A$ with the approximate accuracy $\sigma_{k+1}$ and that $B_k$
does not become ill-conditioned  before $k\leq k_0$. As a result,
the regularized solutions $x^{(k)}$ become better approximations
to $x_{true}$ until iteration $k_0$, and they are deteriorated after that
iteration. At iteration $k_0$, $x^{(k_0)}$ only captures the $k_0$ dominant SVD
components of $A$ and suppress the other $(n-k_0)$ SVD components,
so that it is a best possible regularized solution. As a result,
the pure LSQR has the full regularization for severely
ill-posed problems. We will give a more direct justification on these
assertions in Section 3.3.

\begin{figure}
\begin{minipage}{0.48\linewidth}
  \centerline{\includegraphics[width=7.0cm,height=5cm]{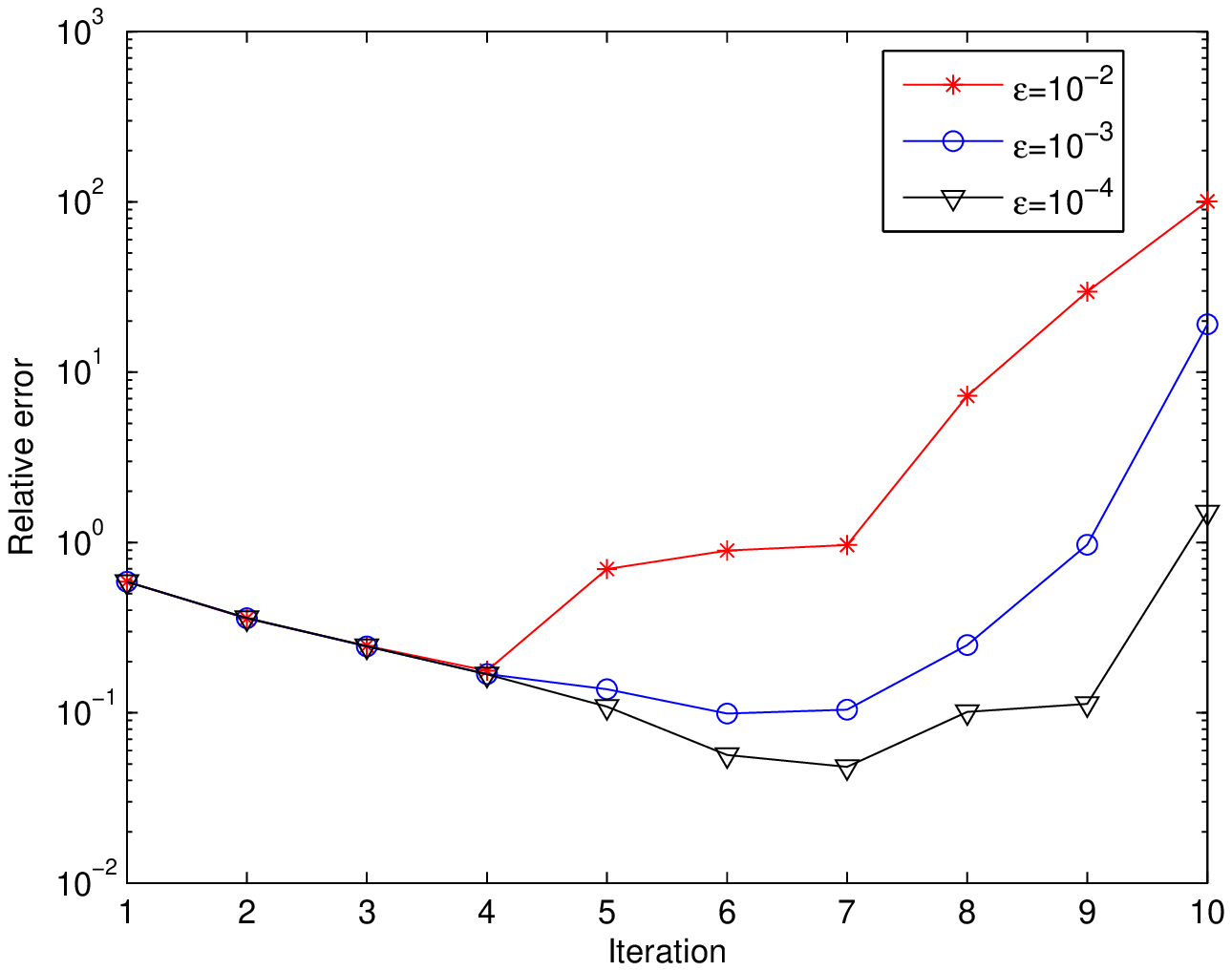}}
  \centerline{(a)}
\end{minipage}
\hfill
\begin{minipage}{0.48\linewidth}
  \centerline{\includegraphics[width=7.0cm,height=5cm]{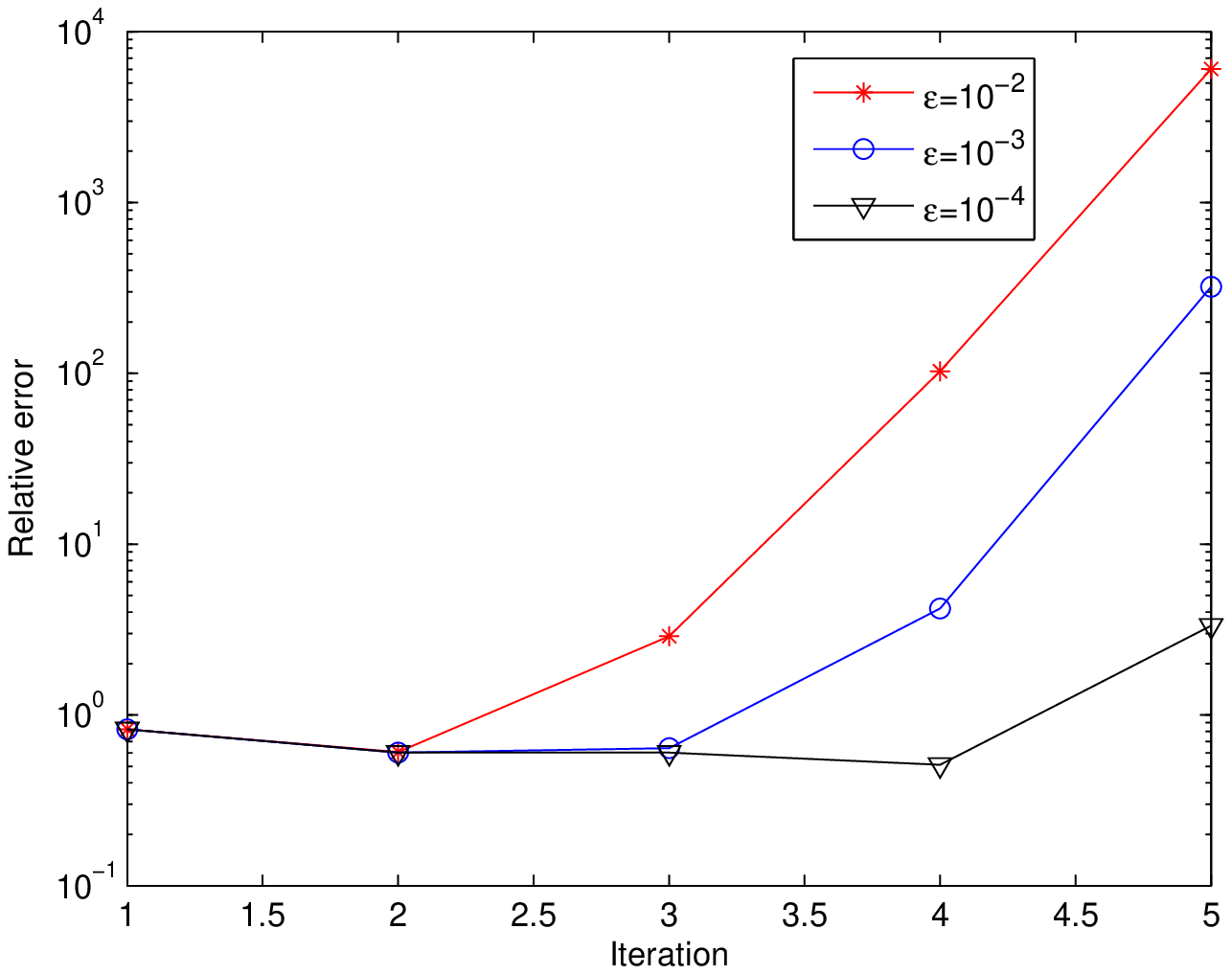}}
  \centerline{(b)}
\end{minipage}
\caption{The relative errors $\left\|x^{(k)}-x_{true}\right\|/
\|x_{true}\|$ with respect to $\varepsilon=10^{-2}, 10^{-3}, 10^{-4}$
for the problems Shaw (left) and Wing (right).}
\label{fig2}
\end{figure}

In Figure~\ref{fig2}, we plot the relative errors
$\left\|x^{(k)}-x_{true}\right\|/\|x_{true}\|$
with different noise levels for these two problems. Obviously,
LSQR exhibits semi-convergence phenomenon. Moreover,
for smaller noise level, we get better regularized solutions
at the cost of more iterations, as expected.

\subsection{Moderately ill-posed problems}

We now consider the following two moderately ill-posed problems \cite{hansen07}.

{\bf Example 3}\ \ \ This problem 'Heat' arises from the inverse heat equation,
and can be obtained by discretizing Volterra integral equation of the first kind,
a class of equations that is
moderately ill-posed, with $[0, 1]$ as integration interval.
The kernel $K(s,t)=k(s-t)$ with
\begin{equation*}\label{}
  k(t)=\frac{t^{-3/2}}{2\sqrt{\pi}}\exp\left(-\frac{1}{4t}\right).
\end{equation*}

{\bf Example 4}\ \ \ This problem is the famous Phillips' test problem. It can be
obtained by discretizing the first kind Fredholm integral
equation \eqref{eq2} with $[-6, 6]$ as both integration and domain intervals.
The kernel $K(s,t)$, the solution $x(t)$ and
the right-hand side $b(s)$ are given by
\begin{equation*}\label{}
  K(s,t)=\left\{\begin{array}{ll} 1+\cos\left(\frac{\pi(s-t)}{3}\right),
  \ \ \ &|s-t|<3;\\ 0,\ \ \ &|s-t|\geq 3, \end{array}\right.
\end{equation*}
\begin{equation*}\label{}
  x(t)=\left\{\begin{array}{ll} 1+\cos\left(\frac{\pi t}{3}\right),
  \ \ \ &|t|<3;\\ 0,\ \ \ &|t|\geq 3, \end{array}\right.
\end{equation*}
\begin{equation*}\label{}
  b(s)=(6-|s|)\left(1+\frac{1}{2}\cos\left(\frac{\pi s}{3}\right)\right)+
  \frac{9}{2\pi}\sin\left(\frac{\pi|s|}{3}\right).
\end{equation*}

 \begin{figure}
\begin{minipage}{0.48\linewidth}
  \centerline{\includegraphics[width=7.0cm,height=5cm]{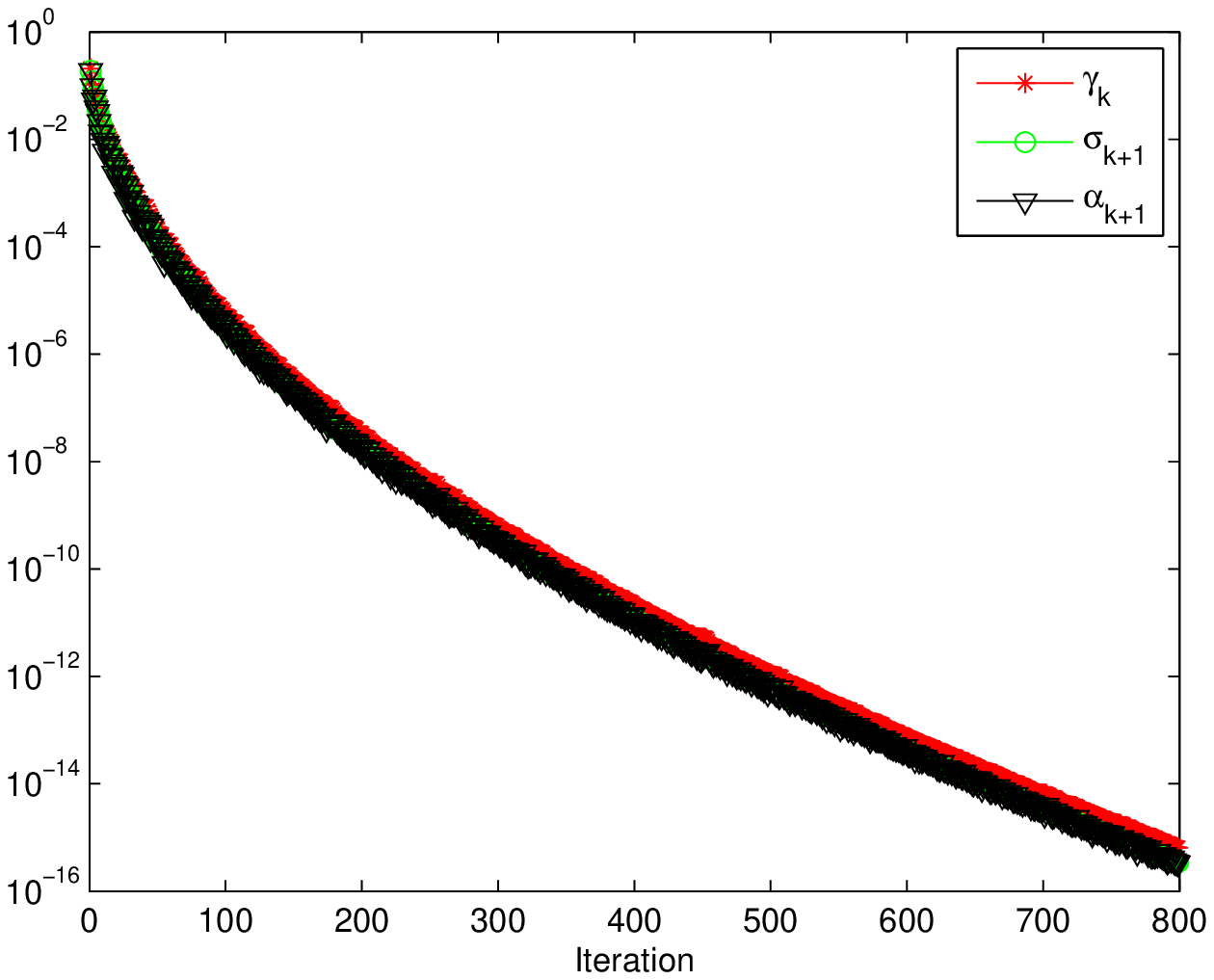}}
  \centerline{(a)}
\end{minipage}
\hfill
\begin{minipage}{0.48\linewidth}
  \centerline{\includegraphics[width=7.0cm,height=5cm]{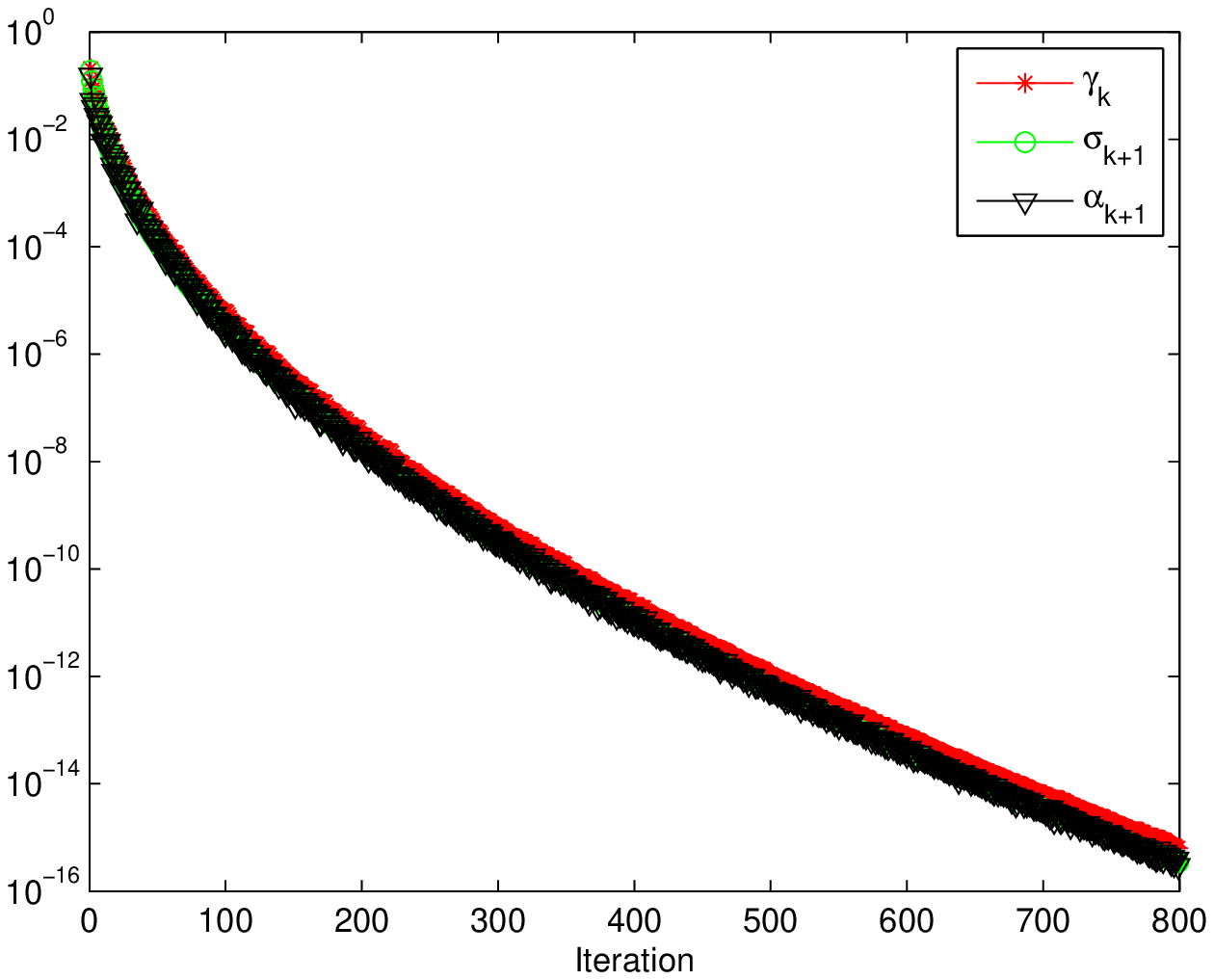}}
  \centerline{(b)}
\end{minipage}
\vfill
\begin{minipage}{0.48\linewidth}
  \centerline{\includegraphics[width=7.0cm,height=5cm]{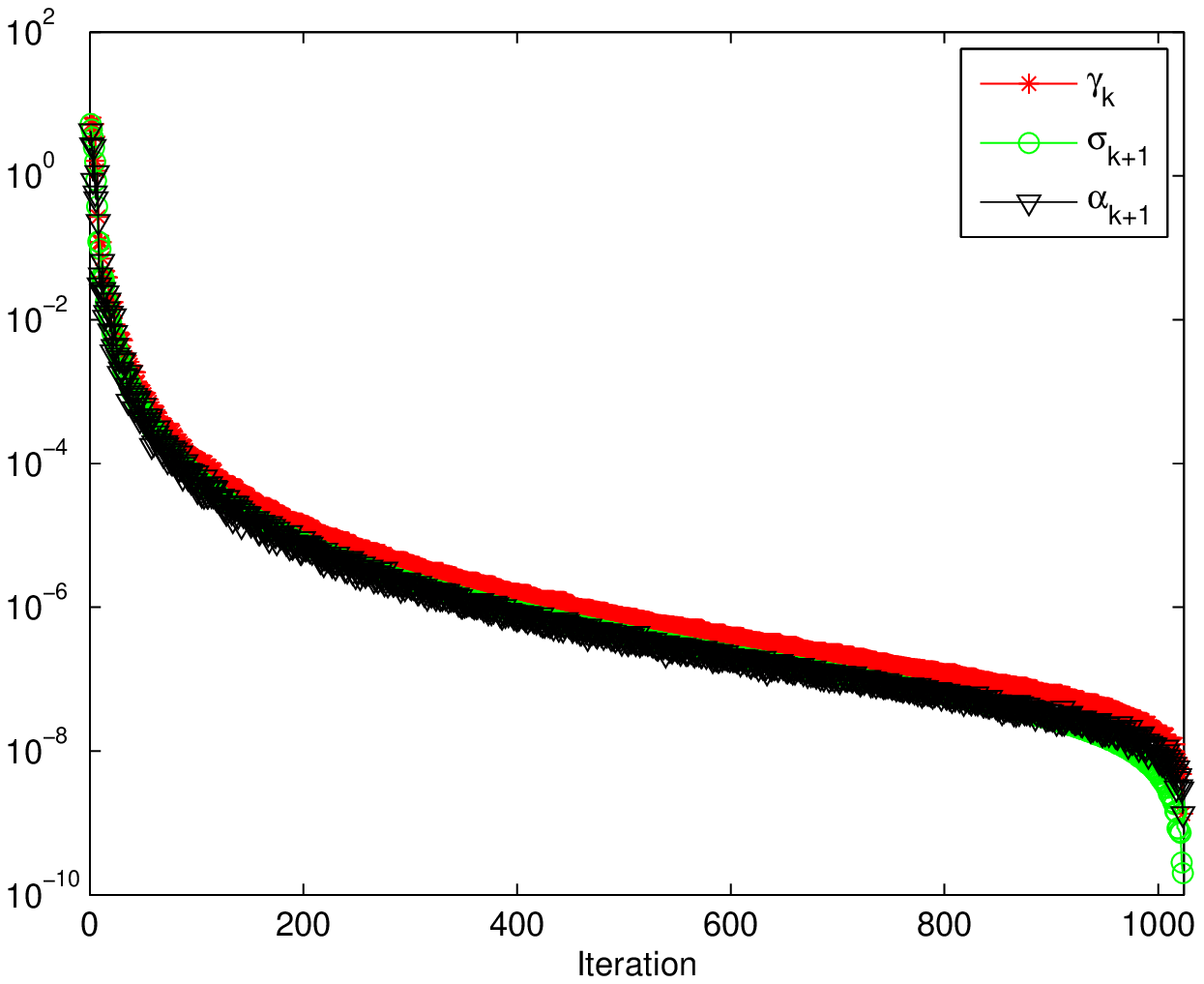}}
  \centerline{(c)}
\end{minipage}
\hfill
\begin{minipage}{0.48\linewidth}
  \centerline{\includegraphics[width=7.0cm,height=5cm]{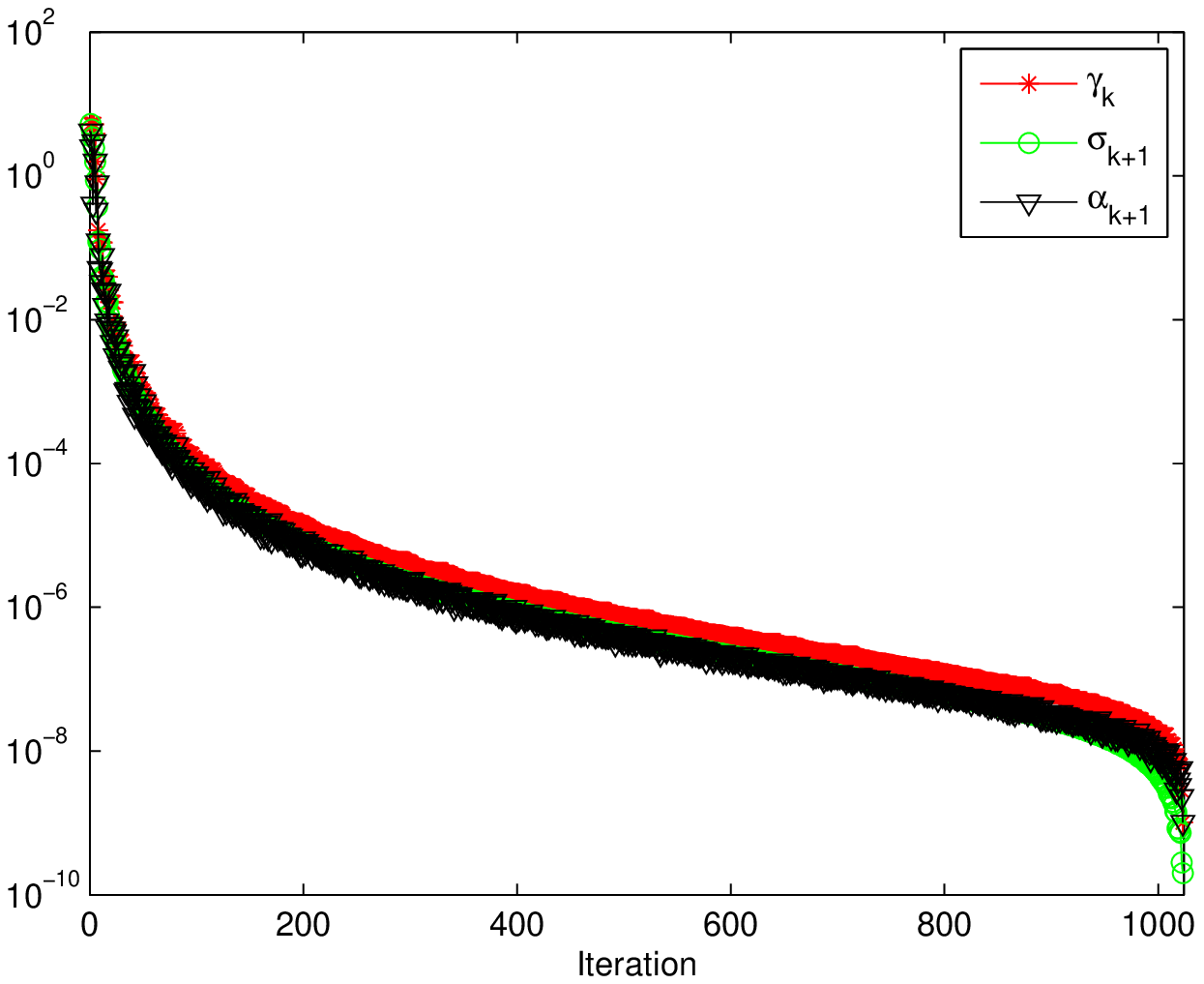}}
  \centerline{(d)}
\end{minipage}
\caption{(a)-(b): Plots of decaying behavior of the sequences $\gamma_k$, $\sigma_{k+1}$ and
$\alpha_{k+1}$ for the problem Heat with $\varepsilon=10^{-2}$ (left)
and $\varepsilon=10^{-3}$ (right); (c)-(d): Plots of decaying behavior of
the sequences $\gamma_k$ and $\sigma_{k+1}$ for the problem Phillips
with $\varepsilon=10^{-3}$ (left) and $\varepsilon=10^{-4}$ (right).}
\label{fig3}
\end{figure}

From Figure~\ref{fig3}, we see that $\gamma_k$ decreases
as fast as $\sigma_{k+1}$, and $\alpha_{k+1}$ decays as fast as $\gamma_k$.
However, slightly different from
severely ill-posed problems, we can observe that the $\gamma_k$ may not be so
close to the $\sigma_{k+1}$, as reflected by the thick rope formed
by three lines. By comparing the behavior of $\gamma_k$
for severely and moderately ill-posed problem, we come to the conclusion
that the $k$-step Lanczos bidiagonalization may generate more accurate
rank $k$ approximation for severely ill-posed problems than for
moderately ill-posed problems, namely, the rank $k$ approximation
$P_{k+1}B_kQ_k^T$ may be more accurate for severely ill-posed problems than
for moderately ill-posed problems. Nonetheless, we have seen
that, for the test moderately ill-posed problems, all the $\gamma_k$ are
still excellent approximations to the $\sigma_{k+1}$, so that
LSQR still has the full regularization.

\begin{figure}
\begin{minipage}{0.48\linewidth}
  \centerline{\includegraphics[width=7.0cm,height=5cm]{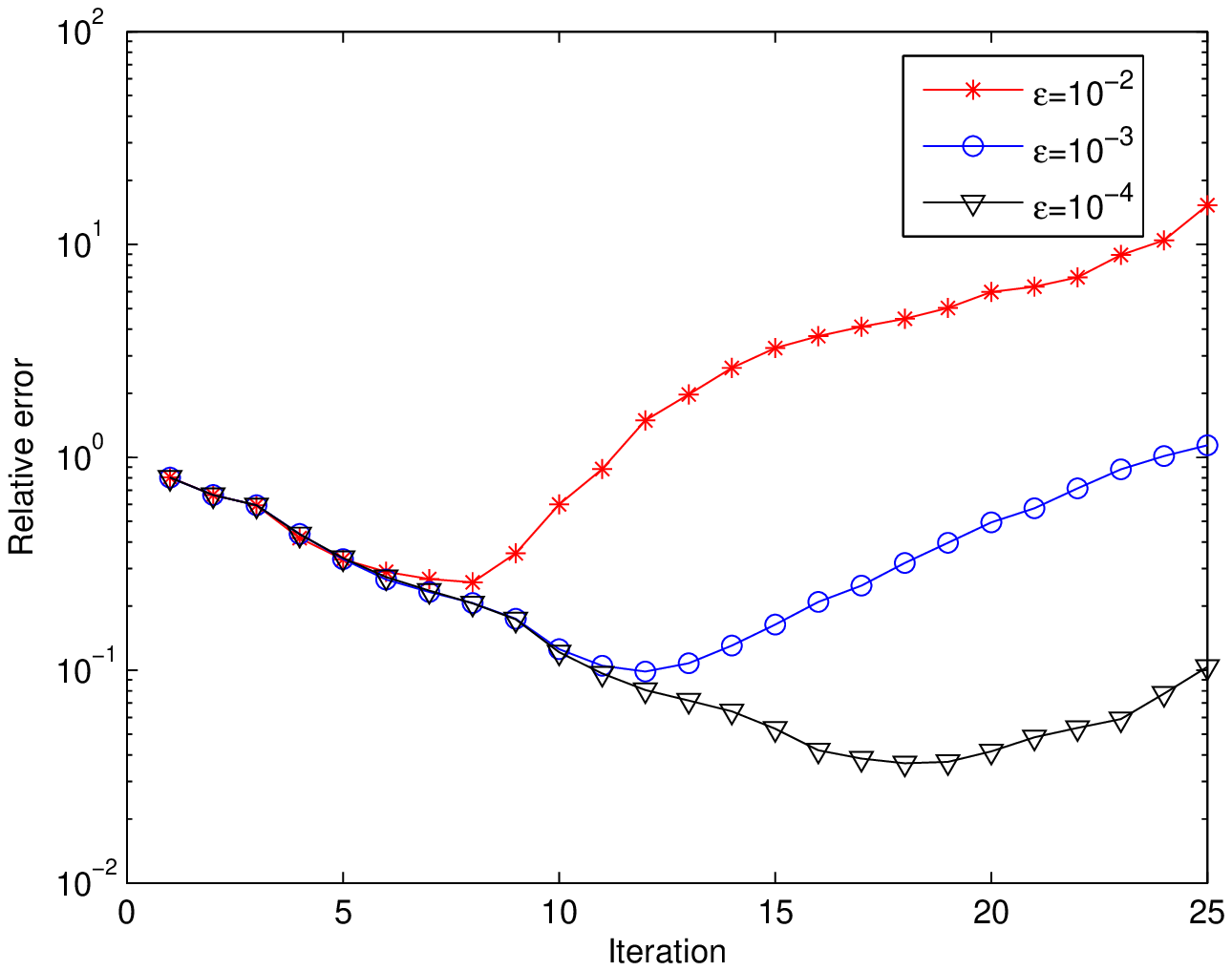}}
  \centerline{(a)}
\end{minipage}
\hfill
\begin{minipage}{0.48\linewidth}
  \centerline{\includegraphics[width=7.0cm,height=5cm]{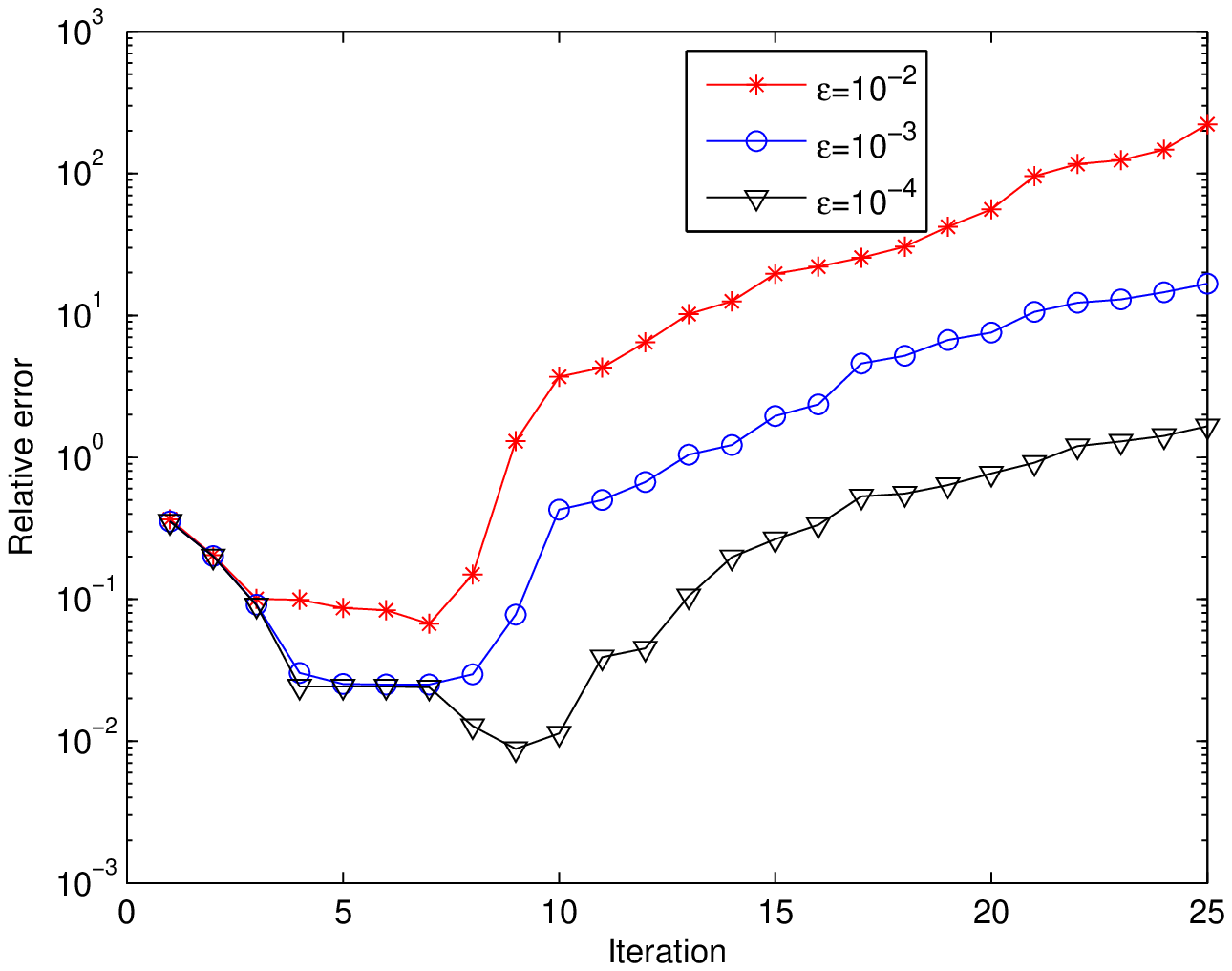}}
  \centerline{(b)}
\end{minipage}
\caption{ The relative errors $\left\|x^{(k)}-x_{true}\right\|/\|x_{true}\|$
with respect to $\varepsilon=10^{-2}, 10^{-3}, 10^{-4}$ for the
problems Heat (left) and Phillips (right).}
\label{fig4}
\end{figure}

In Figure~\ref{fig4}, we depict the relative errors of $x^{(k)}$, and observe
analogous phenomena to those for severely ill-posed problems. A distinction is
that now LSQR needs more iterations for moderately ill-posed problems
with the same noise level.

\subsection{Comparison of LSQR with and without additional TSVD regularization}

For the previous four severely and moderately ill-posed problems,
we now compare the regularizing effects of the pure LSQR
and the hybrid LSQR with the additional TSVD regularization used within
projected problems. We show that LSQR has the full regularization
and no additional regularization is needed, which is based
on the observation that at
semi-convergence the regularized solution by LSQR is as accurate as
that obtained by the hybrid LSQR for each problem.

In the sequel, we only report the results for the noise level
$\varepsilon=10^{-3}$. Results for other $\varepsilon$ are analogous
and thus omitted.

\begin{figure}
\begin{minipage}{0.48\linewidth}
  \centerline{\includegraphics[width=7.0cm,height=5cm]{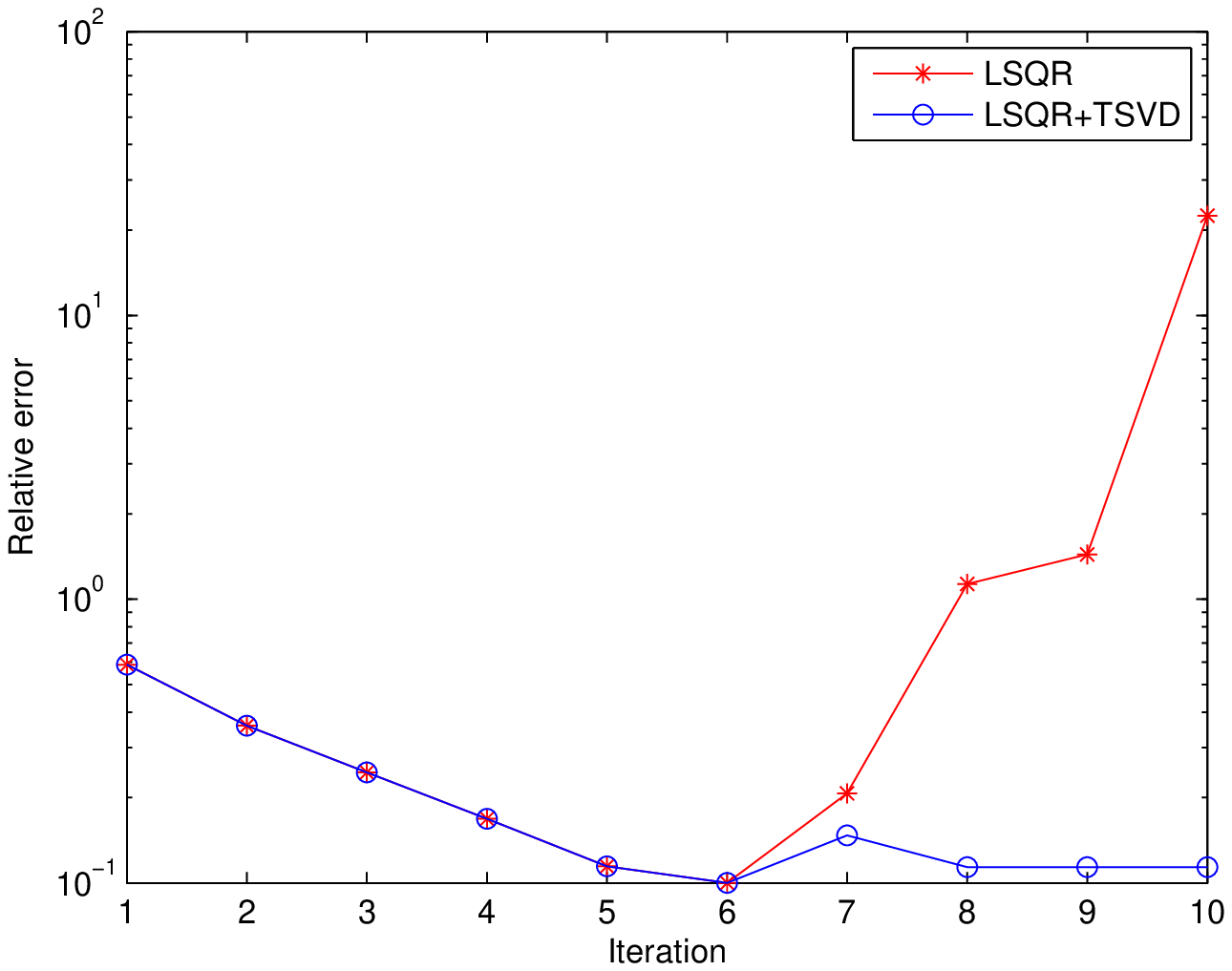}}
  \centerline{(a)}
\end{minipage}
\hfill
\begin{minipage}{0.48\linewidth}
  \centerline{\includegraphics[width=7.0cm,height=5cm]{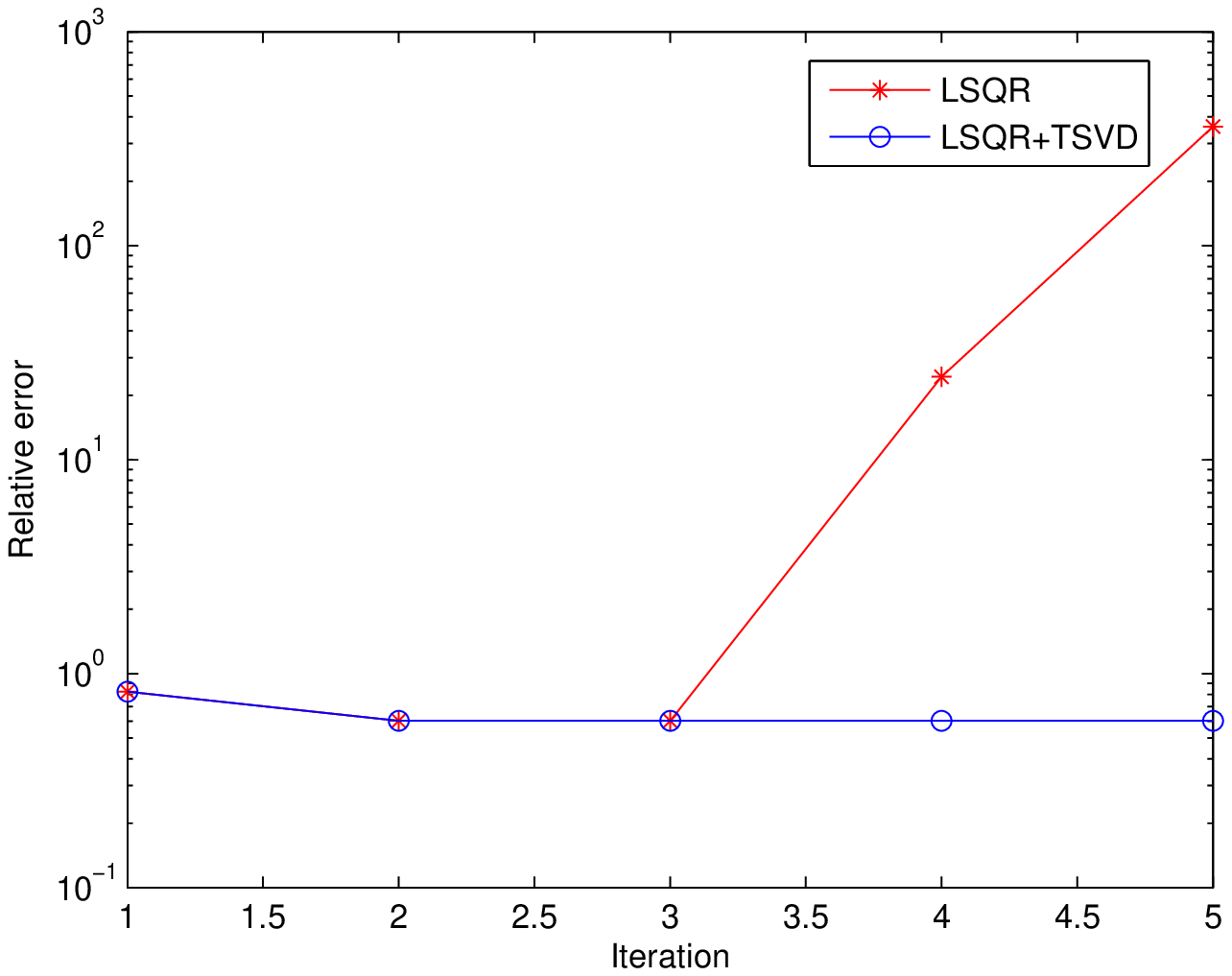}}
  \centerline{(b)}
\end{minipage}
\vfill
\begin{minipage}{0.48\linewidth}
  \centerline{\includegraphics[width=7.0cm,height=5cm]{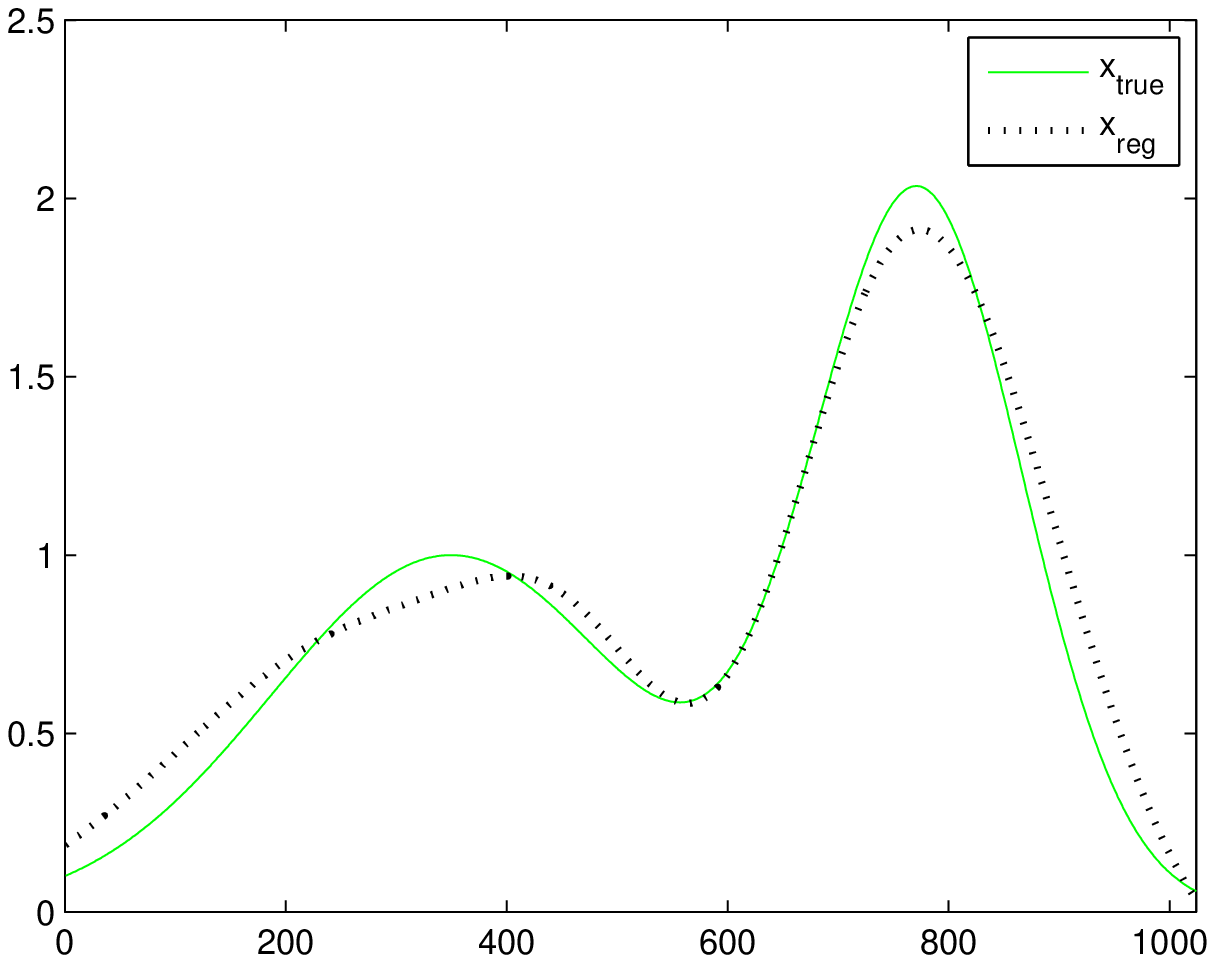}}
  \centerline{(c)}
\end{minipage}
\hfill
\begin{minipage}{0.48\linewidth}
  \centerline{\includegraphics[width=7.0cm,height=5cm]{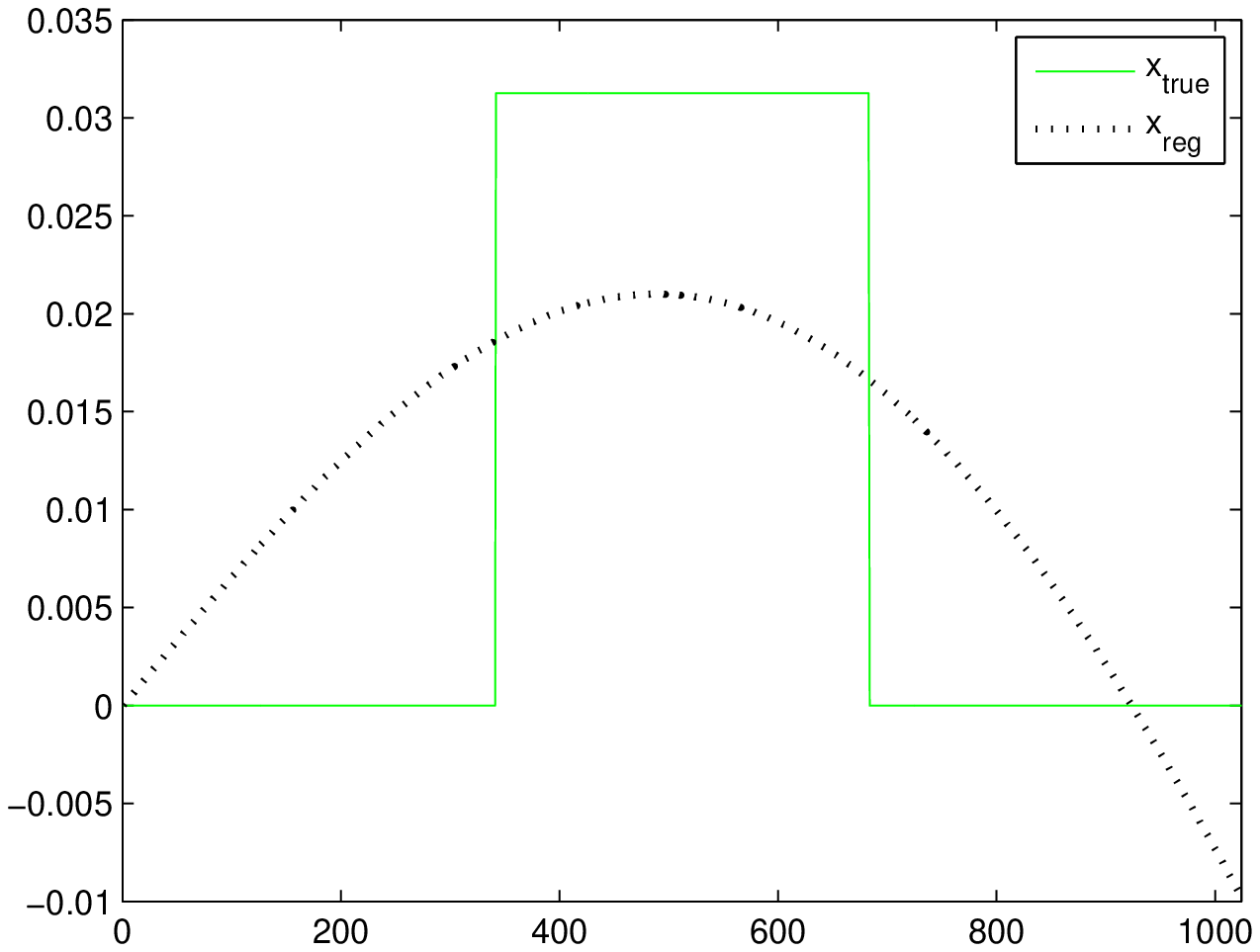}}
  \centerline{(d)}
\end{minipage}
\caption{(a)-(b): The relative errors $\left\|x^{(k)}-x_{true}\right\|/\|x_{true}\|$
with respect to LSQR and LSQR with additional TSVD regularization for
$\varepsilon=10^{-3}$; (c)-(d):
The regularized solutions $x_{reg}$ for the pure LSQR
for the problems Shaw (left) and Wing (right).}
\label{fig5}
\end{figure}

\begin{figure}
\begin{minipage}{0.48\linewidth}
  \centerline{\includegraphics[width=7.0cm,height=5cm]{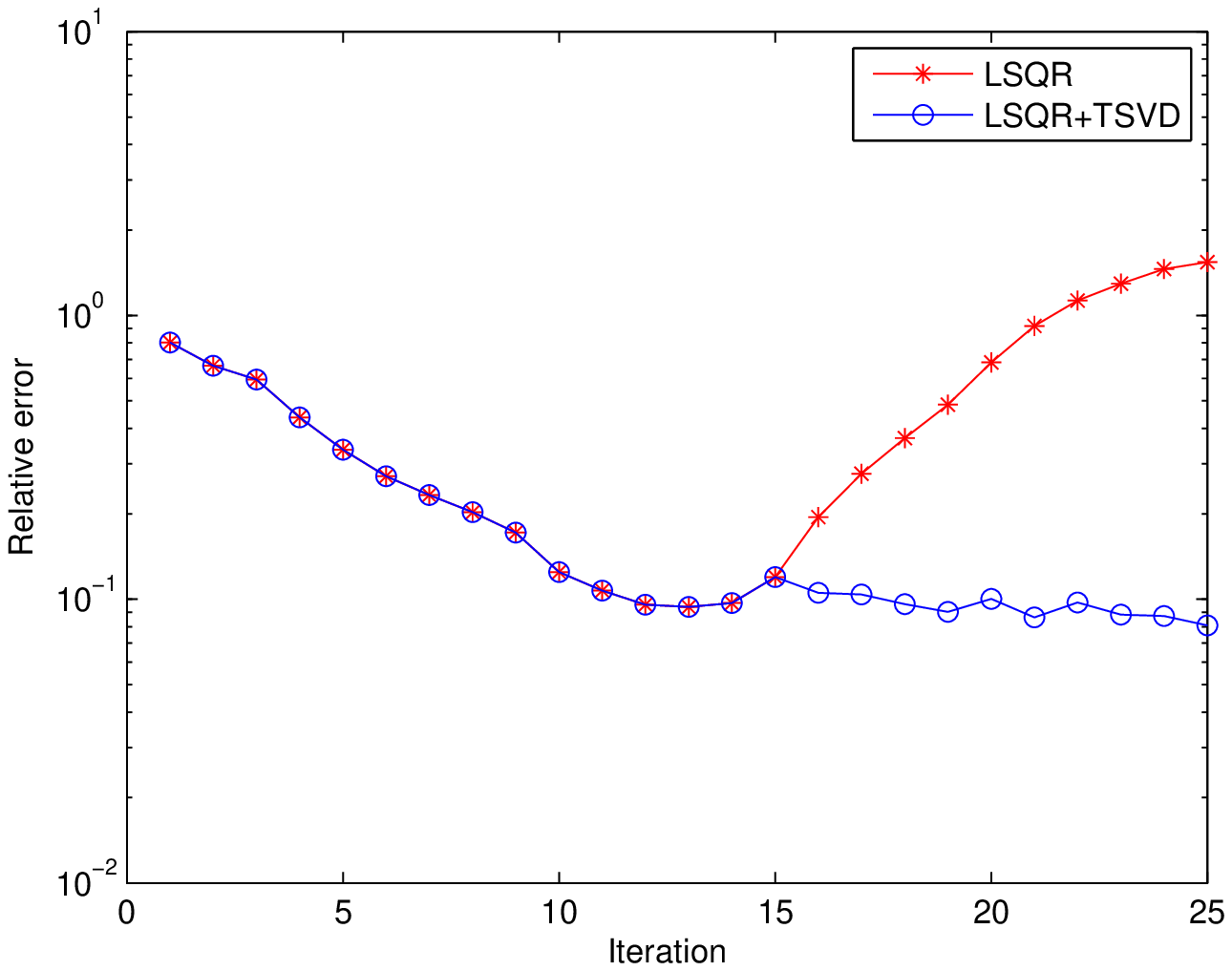}}
  \centerline{(a)}
\end{minipage}
\hfill
\begin{minipage}{0.48\linewidth}
  \centerline{\includegraphics[width=7.0cm,height=5cm]{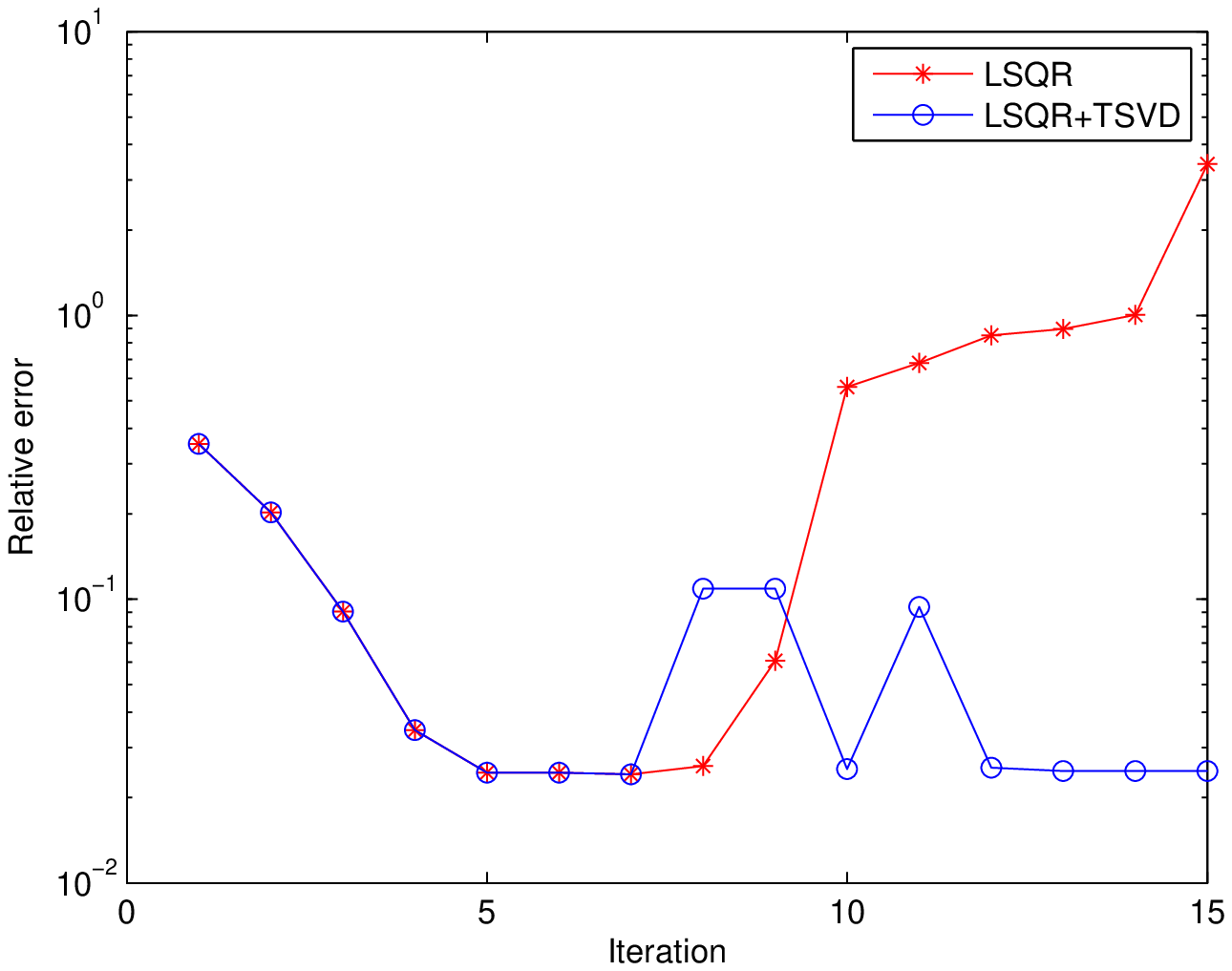}}
  \centerline{(b)}
\end{minipage}
\vfill
\begin{minipage}{0.48\linewidth}
  \centerline{\includegraphics[width=7.0cm,height=5cm]{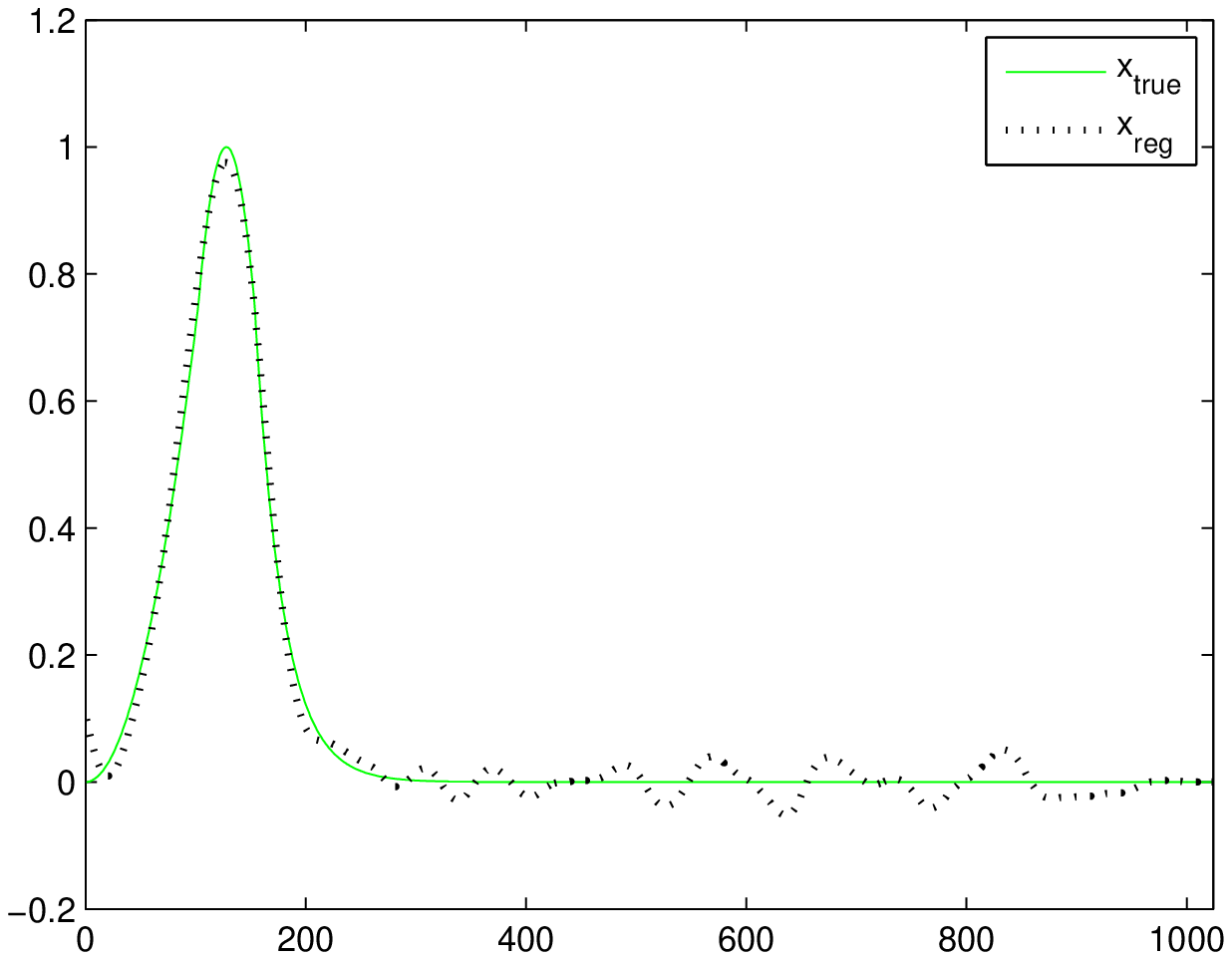}}
  \centerline{(c)}
\end{minipage}
\hfill
\begin{minipage}{0.48\linewidth}
  \centerline{\includegraphics[width=7.0cm,height=5cm]{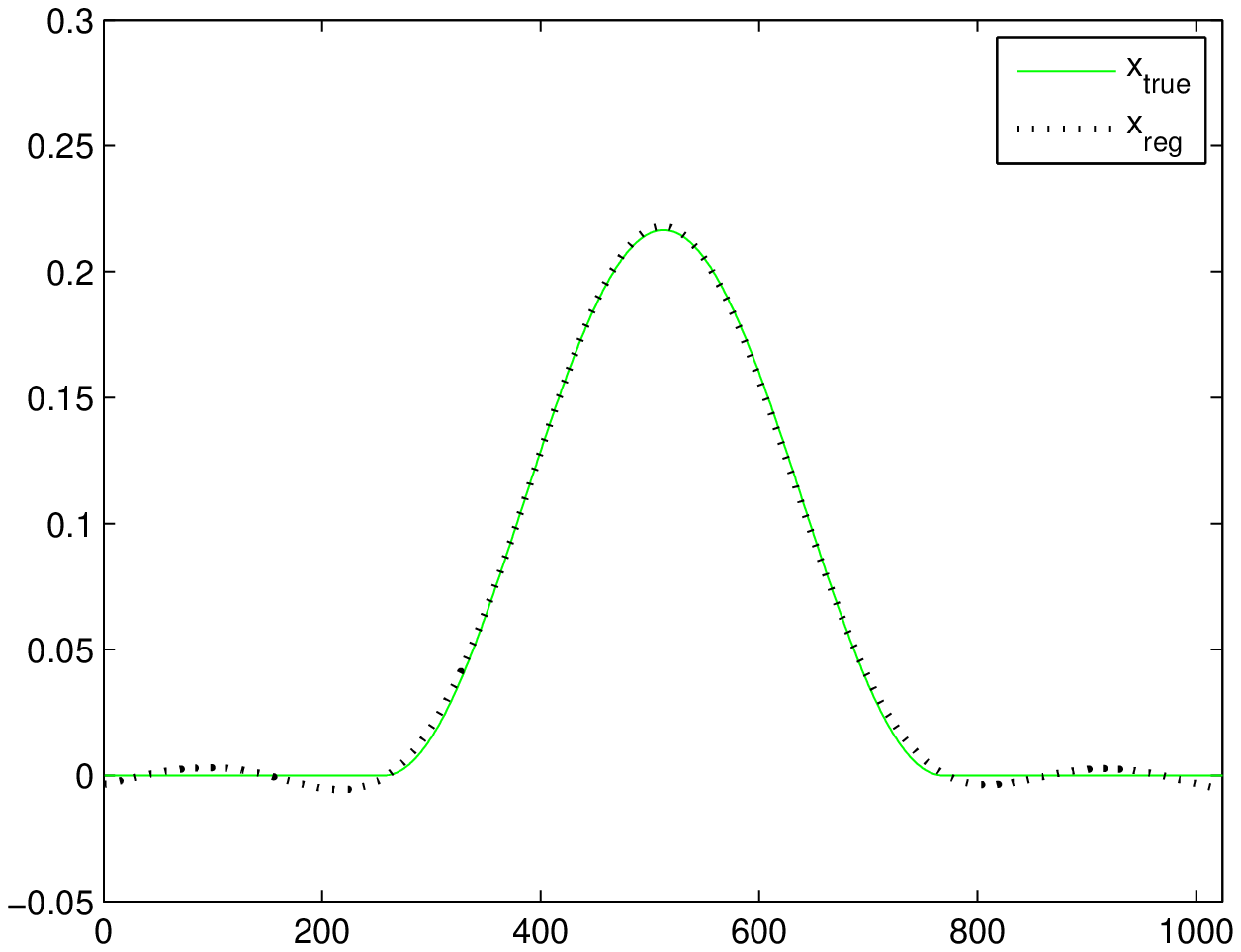}}
  \centerline{(d)}
\end{minipage}
\caption{(a)-(b): The relative errors $\left\|x^{(k)}-x_{true}\right\|/\|x_{true}\|$
obtained by the pure LSQR and LSQR with the additional TSVD regularization for
$\varepsilon=10^{-3}$; (c)-(d):
The regularized solutions $x_{reg}$ for the pure LSQR
for the problems Heat (left) and Phillips (right).}
\label{fig6}
\end{figure}

Figures~\ref{fig5} (a)-(b) and Figures~\ref{fig6} (a)-(b) indicate that the
relative errors of approximate solutions obtained by the two methods reach
the same minimum level, and the hybrid LSQR simply stabilizes the
regularized solutions with the minimum error. This means that
the pure LSQR itself has already found a best possible regularized solution
at semi-convergence and no additional regularization is needed. So it
has the full regularization. Our task is to determine such $k$,
which is the iteration where $\left\|x^{(k+1)}\right\|$ starts to increase
dramatically while its residual norm remains almost unchanged.
The L-curve criterion fits nicely into this task. In these examples, we also
choose $x_{reg}=\arg\min _{k}\left\|x^{(k)}-x_{true}\right\|$ for the pure LSQR.
Figure~\ref{fig5} (c) and Figures~\ref{fig6} (c)-(d) show that the regularized
solutions are generally very good approximations to the true solutions. However,
we should point out that for the problem 'Wing' with a
discontinuous solution, the large relative
error indicates that the regularized solution is a poor approximation
to the true solution, as depicted in Figure~\ref{fig5} (d). Such phenomenon
is due to the fact that the regularization of LSQR and its hybrid variants is
unsuitable for the ill-posed problems with discontinuous solutions. For such kind
of problems, more reasonable regularization is Total Variation Regularization, which
takes the form $\min_{x\in \mathbb{R}^{n}}{\|Ax-b\|^2+\lambda^2\|Lx\|_1^2}$
with $L\not=I$ some $p\times n$ matrix and $\|\cdot\|_1$ the 1-norm \cite{hansen10}.

In what follows, we compare the regularizing effects
of the pure LSQR and hybrid LSQR for mildly ill-posed problems,
showing that LSQR has only the partial
regularization and a hybrid LSQR should be used for this kind of problem
to improve the regularized solution by LSQR at semi-convergence.

{\bf Example 5}\ \ \ The problem 'deriv2' is mildly ill-posed, which is obtained by discretizing
the first kind Fredholm integral equation \eqref{eq2}
with $[0, 1]$ as both integration and domain intervals. The kernel $K(s,t)$ is
Green's function for the second derivative:
\begin{equation*}\label{}
  K(s,t)=\left\{\begin{array}{ll} s(t-1),\ \ \ &s<t;\\ t(s-1),\ \ \
  &s\geq t, \end{array}\right.
\end{equation*}
and the solution $x(t)$ and
the right-hand side $b(s)$ are given by
\begin{equation*}\label{}
  x(t)=\left\{\begin{array}{ll} t,
  \ \ \ &t<\frac{1}{2};\\ 1-t,\ \ \ &t\geq\frac{1}{2}, \end{array}\right. \ \ \
  b(s)=\left\{\begin{array}{ll} (4s^3-3s)/24,\ \ \ &s<\frac{1}{2};
  \\ (-4s^3+12s^2-9s+1)/24,\ \ \ &s\geq\frac{1}{2}. \end{array}\right.
\end{equation*}

\begin{figure}
\begin{minipage}{0.48\linewidth}
  \centerline{\includegraphics[width=7.0cm,height=5cm]{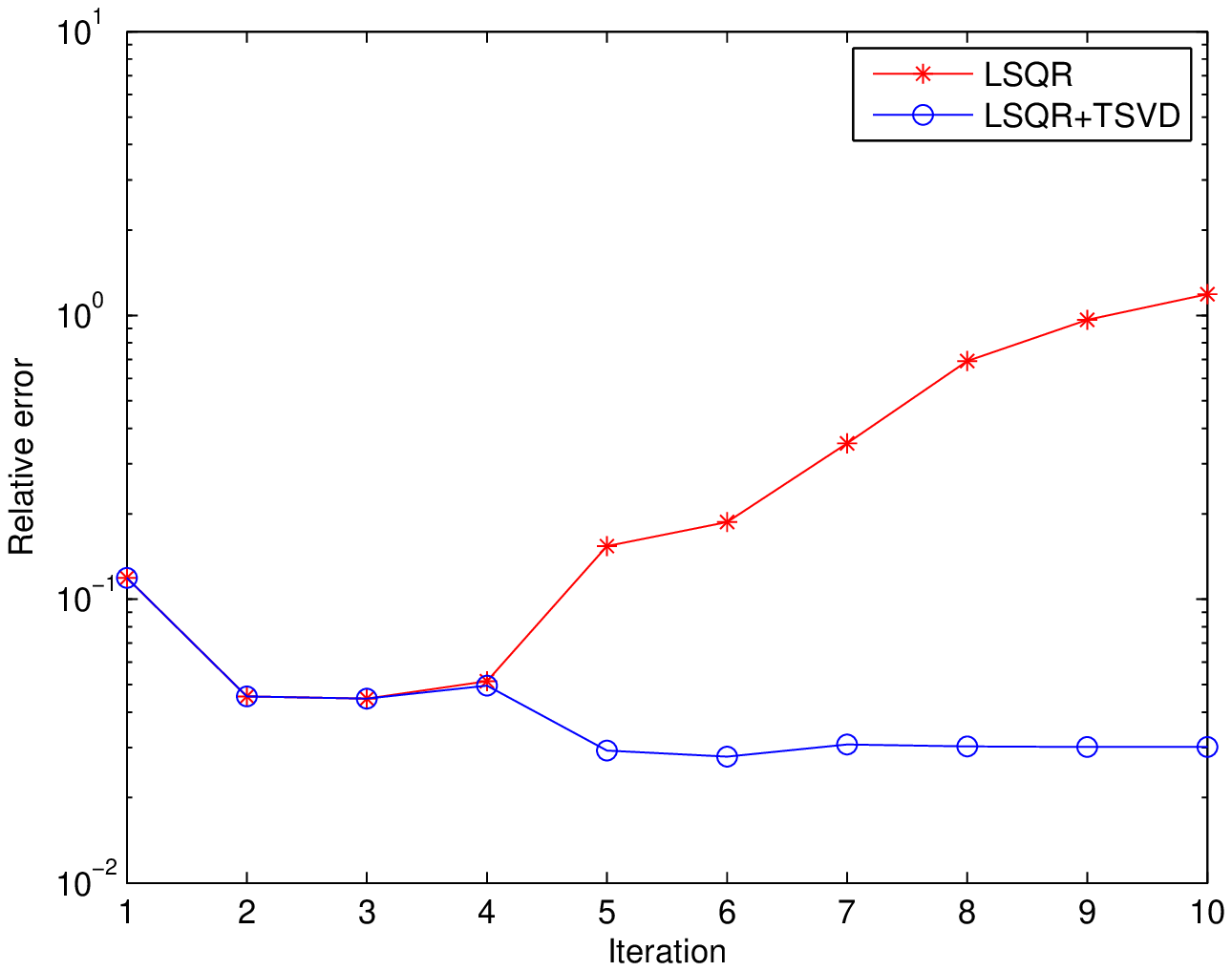}}
  \centerline{(a)}
\end{minipage}
\hfill
\begin{minipage}{0.48\linewidth}
  \centerline{\includegraphics[width=7.0cm,height=5cm]{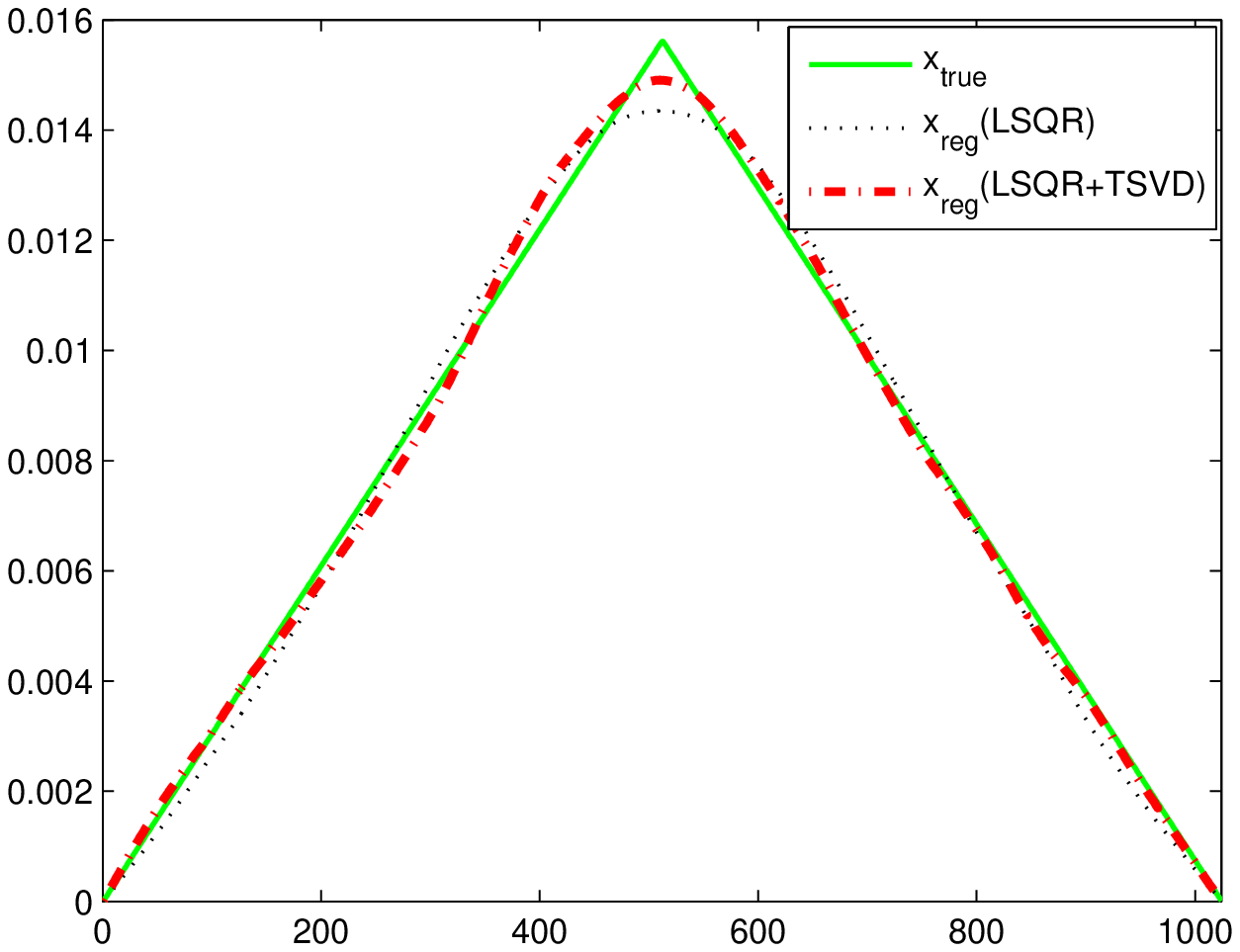}}
  \centerline{(b)}
\end{minipage}
\caption{ The relative errors $\left\|x^{(k)}-x_{true}\right\|/\|x_{true}\|$ and
the regularized solution $x_{reg}$ with respect to LSQR and
LSQR with the additional TSVD regularization for the problem Deriv2 and
$\varepsilon=10^{-3}$.}
\label{fig7}
\end{figure}

Figure~\ref{fig7} (a) shows that the relative errors
of approximate solutions by the hybrid LSQR reach a considerably
smaller minimum level than those by the pure LSQR, a clear indication
that LSQR has the partial regularization. As we have seen, the hybrid LSQR
expands the Krylov subspace until it contains enough
dominant SVD components and, meanwhile, additional regularization effectively
dampen the SVD components corresponding
to small singular values. For instance, the semi-convergence of the pure LSQR
occurs at iteration $k=3$, but it is not enough. As the hybrid LSQR
shows, we need a larger six dimensional Krylov subspace $\mathcal{K}_6(A^TA,A^Tb)$
to construct a best possible regularized solution. We also
choose $x_{reg}=\arg\min _{k}\left\|x^{(k)}-x_{true}\right\|$ for the pure LSQR
and the hybrid LSQR. Figure~\ref{fig7} (b) indicates that the regularized
solution obtained by the hybrid LSQR is a considerably better approximation
to $x_{true}$ than that by the pure LSQR, especially in the non-smooth middle
part of $x_{true}$.

\section{Conclusions}\label{SectionCon}

For large-scale discrete ill-posed problems, LSQR and CGLS are commonly
used methods. These methods have regularizing effects and exhibit
semi-convergence. However, if a small Ritz value appears before the
methods capture all the needed dominant
SVD components, the methods have only the partial regularization and
must be equipped with additional regularization so that best possible
regularized solutions can be found. Otherwise, LSQR has the full
regularization and can compute best possible regularized solutions without
additional regularization needed.

We have proved that the underlying $k$-dimensional Krylov subspace captures
the $k$ dimensional dominant right singular space better for
severely and moderately ill-posed problems than for mildly ill-posed problems.
This makes LSQR have better regularization for the first two kinds of problems
than for the third kind. Furthermore,  we have shown that LSQR
generally has only the partial regularization
for mildly ill-posed problems. Numerical experiments have demonstrated
that LSQR has the full regularization for severely and moderately ill-posed
problems, stronger than our theory predicts, and it has the
partial regularization for mildly moderately ill-posed problems,
compatible with our assertion. Together with
the observations \cite{bazan14,gazzola14,gazzola-online},
it appears that the excellent performances of LSQR
on severely and moderately ill-posed problems generally hold.

As for future work, it is more appealing to derive an accurate estimate for
$\|\Delta_k\|$ other than $\|\Delta_k\|_F$,
as it plays a crucial role in analyzing the accuracy $\gamma_k$ of the
rank $k$ approximation, generated by Lanczos bidiagonalization, to $A$.
Accurate bounds for $\gamma_k$ are
the core of completely understanding the regularizing effects of LSQR, but
our bound \eqref{gammak} for $\gamma_k$ is conservative and is expected to
be improved on substantially.
Since CGLS is mathematically equivalent to LSQR, our results apply to CGLS
as well. Our current work has helped to better understand the regularization
of LSQR and CGLS. But for a complete understanding of the
intrinsic regularizing effects of LSQR and CGLS, we still have
a long way to go, and more research is needed.

\section*{Acknowledgements} We thank the three referees very much for their
valuable suggestions and comments, which made us improve the presentation
of the paper.

\end{document}